\documentclass[11pt, reqno]{amsart}
\usepackage{lmodern}
\usepackage{amsmath, amsthm, amssymb, amsfonts}
\usepackage[normalem]{ulem}
\usepackage{hyperref}
\usepackage[all,cmtip]{xy}
\usepackage{verbatim}
\usepackage{nccmath}
\usepackage{stmaryrd}
\usepackage{caption}
\usepackage{mathrsfs}
\usepackage{mathtools}
\usepackage{esvect}
\usepackage{cite}
\usepackage{bbm}
\usepackage{eucal}

\usepackage{mathbbol}
\usepackage{tabularx}
\usepackage[toc,page]{appendix}
\makeatletter
\makeatletter

\usepackage{tikz-cd}

\theoremstyle{plain}
\newtheorem{thm}{Theorem}[section]
\newtheorem{cor}[thm]{Corollary}

\newtheorem{lemma}[thm]{Lemma}
\newtheorem{prop}[thm]{Proposition}

\newtheorem{claim}{Claim}

\newtheorem{conjecture}{Conjecture}  

\newtheorem{thml}{Theorem}

\theoremstyle{definition}
\newtheorem{defn}[thm]{Definition}

\makeatletter
\newcommand\ackname{Acknowledgements}
\if@titlepage
\newenvironment{acknowledgements}{%
	\titlepage
	\null\vfil
	\@beginparpenalty\@lowpenalty
	\begin{center}%
		\bfseries \ackname
		\@endparpenalty\@M
\end{center}}%
{\par\vfil\null\endtitlepage}
\else

\fi
\makeatother

\theoremstyle{remark}
\newtheorem{rmk}[thm]{Remark}

\newcommand{\BA}{{\mathbb{A}}}

\newcommand{\BC}{{\mathbb{C}}}

\newcommand{\BQ}{{\mathbb{Q}}}
\newcommand{\BR}{{\mathbb{R}}}

\newcommand{\BT}{{\mathbb{T}}}

\newcommand{\BZ}{{\mathbb{Z}}}

\newcommand{\CA}{{\mathcal A}}
\newcommand{\CB}{{\mathcal B}}

\newcommand{\CE}{{\mathcal E}}
\newcommand{\CF}{{\mathcal F}}
\newcommand{\CG}{{\mathcal G}}
\newcommand{\CH}{{\mathcal H}}

\newcommand{\CK}{{\mathcal K}}
\newcommand{\CL}{{\mathcal L}}

\newcommand{\CN}{{\mathcal N}}
\newcommand{\CO}{{\mathcal O}}

\newcommand{\CQ}{{\mathcal Q}}

\newcommand{\FM}{{\mathfrak{M}}}
\newcommand{\FN}{{\mathfrak{N}}}

\newcommand{\vsf}{{\mathsf{v}}}
\newcommand{\wsf}{{\mathsf{w}}}

\newcommand{\rsf}{{\mathsf{r}}}
\newcommand{\dsf}{{\mathsf{d}}}
\newcommand{\usf}{{\mathsf{u}}}
\newcommand{\asf}{{\mathsf{a}}}
\newcommand{\msf}{{\mathsf{m}}}

\newcommand{\ksf}{{\mathsf{k}}}

\newcommand{\Usf}{{\mathsf{U}}}

\newcommand{\ch}{{\mathrm{ch}}}

\newcommand{\rk}{{\mathrm{rk}}}

\newcommand{\Cohc}{\mathrm{Coh}}

\newcommand{\VW}{\mathsf{VW}}

\newcommand{\QM}{\mathsf{QM}}

\newcommand{\rel}{\mathrm{rel}}

\newcommand{\tr}{{\mathrm{tr}}}

\DeclareFontFamily{OT1}{rsfs}{}
\DeclareFontShape{OT1}{rsfs}{n}{it}{<-> rsfs10}{}
\DeclareMathAlphabet{\curly}{OT1}{rsfs}{n}{it}

\newcommand\Hom{\operatorname{Hom}}
\newcommand{\p}{\mathbb{P}}

\newcommand{\Jac}{{\mathrm{Jac}}}

\newcommand{\td}{\mathrm{td}}

\newcommand{\Pic}{\mathop{\rm Pic}\nolimits}

\newcommand{\GW}{\mathsf{GW}}

\newcommand{\ev}{{\mathrm{ev}}}

\newcommand{\thickslash}{\mathbin{\!\!\pmb{\fatslash}}}

\newcommand{\Ob}{\mathrm{Ob}}
\newcommand{\M}{\mathsf{M}}

\begin{document}

	\title[On quasimap invariants of moduli spaces of Higgs bundles]
	{On quasimap invariants of moduli spaces of Higgs bundles}

	\author{Denis Nesterov}
	\address{University of Vienna, Faculty of Mathematics}
	\email{denis.nesterov@univie.ac.at}
	\maketitle
	\begin{abstract}
		We compute odd-degree genus 1 quasimap (and Gromov--Witten) invariants of moduli spaces of Higgs $\mathrm{SL}_2$-bundles on a curve of genus $g\geq2$.  We also compute certain invariants for all prime ranks. This proves some parts of author's conjectures on quasimap invariants of moduli spaces of Higgs bundles.
		 
		  More generally, our methods provide a computation scheme for genus 1 quasimap (and Gromov--Witten)  invariants in the case when degrees of maps are coprime to the rank. This requires a careful analysis of the localisation formula for certain Quot schemes parametrising higher-rank quotients on an elliptic curve.  Invariants for degrees which are not coprime to the rank exhibit a very different structure for a reason that we explain. 
		\end{abstract}
		\setcounter{tocdepth}{1}
	\tableofcontents
	\section{Introduction}
\subsection{Quasimaps}
 Let $M(\dsf)$ be a moduli space of semistable  Higgs $\mathrm {SL_2}$-bundles of degree $\dsf$ on a curve $C$ of genus $g\geq 2$.  In this work, we consider quasimaps from a fixed elliptic curve $E$ to $M(\dsf)$. These are maps from $E$ to the stack of all Higgs sheaves mapping generically to $M(\dsf)$. 

The (reduced) expected dimension of a moduli space of quasimaps  up to translations by $E$  is 0, hence by \cite{NHiggs} it produces an invariant 
\[ \QM^\bullet_{\dsf,\wsf} \in \BQ,\]
where $\wsf \in \BZ$ is the degree\footnote{The Picard rank of $M(\dsf)$ is 1, the degree is taken with respect to the ample generator of $\Pic(M(\dsf))$. } of quasimaps. 
Assuming $\dsf=1$ (or equivalently $\dsf$ is odd),  we determine these invariants for odd degrees $\wsf$.  Let 
\[ \Usf(q):= \log \left(\prod_{k\geq 0}(1-q^k)\right). \] 
\begin{thml}[Theorem \ref{themain}] \label{theeq2} \
\[
\sum_{\mathrm{odd }\ \wsf }\QM^\bullet_{1,\wsf}q^\wsf= (2-2g)2^{2g-1} \left( \Usf(q)-\Usf(-q)\right).\]
\end{thml}
By \cite[Corollary 10.12]{NHiggs}, this also determines genus 1 Gromov--Witten invariant $\GW^\bullet_{1,\wsf}$ of $M(1)$, since
 \[ \QM^\bullet_{1,\wsf}= \GW^\bullet_{1,\wsf},\]
 if $\wsf$ is odd. The invariant $\GW^\bullet_{1,\wsf}$ is defined analogously but via the moduli space of stable maps from $E$.  Moreover, by Corollary \ref{comp}, these quasimap invariants determine certain Vafa--Witten invariants with insertions on the product  $C\times E$. 
\subsection{Degree 0 Higgs bundles} In \cite{NHiggs}, a notion of extended degree\footnote{The extended degree aims to capture presence of torsion classes in the cohomology of moduli spaces of $\mathrm{PGL}_\rsf$-bundles.  It is essential for the formulation of enumerative mirror symmetry.} $(\wsf,\asf) \in \BZ\oplus \BZ_2$ of quasimaps to $M(\dsf)$ was defined (see Definition \ref{defn}), as well as the associated invariants for an arbitrary $\dsf$,
\[ \QM^{\asf, \bullet}_{\dsf,\wsf} \in \BQ.\]
 If $\dsf=1$, the extended degree is determined by the parity of $\wsf$, i.e.\  $(\wsf,\asf)=(\wsf, [\wsf]_2)$, where $[\wsf]_2:=\wsf \textrm{ mod }2.$ In particular, for odd $\wsf$ we have 
 \[ \QM^{1, \bullet}_{1,\wsf}= \QM^\bullet_{1,\wsf}.\]
In this work, we also determine the quasimap invariants $\QM^{1, \bullet}_{0,\wsf}$ associated to a moduli space of degree 0 Higgs $\mathrm{SL}_2$-bundles $M(0)$.
\begin{thml}[Theorem \ref{themain}] \label{theeq} \
\begin{equation*} 
\sum_\wsf \QM^{1, \bullet}_{0,\wsf} q^\wsf= (2-2g)2^{2g-1} \left( \Usf(q)+\Usf(-q)\right)  .
\end{equation*}
\end{thml}
 As it was argued in \cite{NHiggs}, these invariants can be seen as Gromov--Witten type invariants of the stack $M(0)$. Moreover, we compute certain invariants for $\wsf=0$ in Corollary \ref{const}.
 Theorem \ref{theeq2} and \ref{theeq} with Corollary \ref{const} confirm parts of \cite[Conjecture A, B]{NHiggs},  providing a solid evidence for \cite[Conjecture A, B, C]{NHiggs}.

 We want to draw the reader's attention to a peculiar coincidence. Genus 1 positive-degree Gromov-Witten invariants of an elliptic curve $E$ are given by the following generating series
\[-\Usf(q)=- \log\left( \prod_{k\geq 0}(1-q^k)\right), \]
as is shown in \cite{Di}. Our methods make this coincidence in some sense less surprising\footnote{However, we do not claim that we can fully explain this coincidence, so we urge the reader to treat this sentence as mostly rhetoric.}, because in fact the roles of $C$ and $E$ can be exchanged, allowing us to treat invariants $\QM^{\asf, \bullet}_{\dsf,\wsf}$ in terms of other invariants which are expressible via $E$ alone.   We now explain how this is done (see also Remark \ref{genus1}).
\subsection{Methods} Correspondence between quasimap invariants of $M(\dsf)$ and Vafa--Witten invariants of $C\times E$, discussed in \cite{NHiggs}, is essential. Our argument uses a combination of 
\begin{itemize} \item wall-crossing for Vafa--Witten invariants,
	\item quasimap wall-crossing.
\end{itemize} 
 The wall-crossing for Vafa--Witten invariants is conjectured to be trivial for a complex surface $S$ with $p_g(S)> 0$ or equivalently with $b_2^+(S)> 1$. We will sketch an argument for its triviality, assembling various results from the existing literature. However, these results usually assume that $b_1(S)=0$, mainly in order to simplify the exposition. In our case, $S$ is a product of two non-rational curves, hence $b_1(S)\neq 0$. We therefore make an assumption that existing results extend to the case of $S$ with $b_1(S)\neq 0$ . See Section \ref{secvw} for more details. 
 
   Changing stability on $C\times E$ from the one that has a high degree on $E$ to the one that has a high degree on $C$ corresponds to passing from quasimaps  $E \dashrightarrow M(\dsf)$ to quasimaps $C \dasharrow M'(\asf)$, where  $M'(\asf)$ is a moduli space of Higgs $\mathrm{SL}_2$-bundles on $E$. Since the wall-crossing for Vafa--Witten invariants is trivial,  this gives rise to an equivalence of associated invariants - Corollary \ref{comp1} and \ref{comp}.

In this way, invariants $\QM^{1, \bullet}_{\dsf,\wsf}$ correspond to quasimap invariants of $M'(1)$. This makes computation more accessible, because $M'(1)$ is just a point (there is a unique stable Higgs $\mathrm{SL}_2$-bundle of degree 1 on $E$). Hence the corresponding quasimap invariants can be effectively computed by the quasisection wall-crossing - they will be equal to the wall-crossing invariants of quasisection wall-crossing formula, which are just  Euler characteristics of certain Quot schemes on $E$ (quotiented by the action of $E$). The summary of preceding discussion is depicted in Figure \ref{square}. 
   \begin{figure}[h!] \vspace{0.7cm} 
\[  
   \boxed{\mathbf{\#}^{\mathrm{vir}}\{E \dashrightarrow M(\dsf)\} \approx  \#^{\mathrm{vir}} \{C \dashrightarrow M'(\asf)\} \approx \#^{\mathrm{vir}} \{ G \twoheadrightarrow  Q  \text{ on }E \}}
 \]
   	\caption{Summary}
   \label{square}
   \vspace{-0.5cm}
   \end{figure}

Complications arise due to the fact that in reality one needs to consider quasisections of $C$ to $M'(\asf)$ instead of quasimaps, this is the reason we did not put the sign of equalities above (another reason is that we have to quotient by $E$). Moduli spaces of quasisections and quasimaps are essentially isomorphic in this case, but the obstruction theories are not. This is the main source of technicalities in our calculations, as the obstruction theory of quasisections no longer has surjective cosections (they receive a certain twist). Nevertheless, we obtain the same vanishing results as in the case of an obstruction theory with a surjective cosection, Theorem \ref{thetheorem}. We refer to Section \ref{squasisections} for more on quasisections in our context. Quasisections are treated in greater detail for arbitrary smooth fibrations  in \cite{LeeW}. 

\subsection{Higher rank}

  Almost everything presented in this work applies to an arbitrary rank $\rsf$, except the following two results. Firstly, Claim \ref{claim} is stated for a prime rank $\rsf$, because of Thomas' vanishing result \cite[Corollary 5.30]{Th20}. Secondly, the analysis of the wall-crossing Quot schemes in Section \ref{secquot} is done only for $\rsf=2$.
  
   The case of $\rsf>2$ requires Quot schemes parametrizing  quotients of higher rank on $E$. One can always deform to a sum of line bundle and use torus-localisations. However, since a sum of line bundles is not stable and we consider higher-rank quotients, the resulting Quot schemes have non-trivial obstruction theories, which slightly obscures localisation formulas (at least for the author), hence this will be addressed elsewhere. To this end, we conjecture an expression for higher rank invariants in Conjecture \ref{conj} and provide a basic check of the conjecture, Proposition \ref{themain2}. 
\subsection{Even degrees} There is a good reason  why we cannot compute invariants for even degrees $\wsf$ (or more generally, for degrees coprime to the rank) using the same methods. This case corresponds to the moduli space of  degree 0 Higgs $\mathrm{SL}_2$-bundles $M'(0)$  on $E$. The space $M'(0)$ is no longer a point. In fact, all Higgs bundles of degree 0 on $E$ are strictly semistable and are given by direct sums of degree 0 line bundles. As such, the moduli space $M'(0)$ is not complicated, but since it is a stack, its moduli spaces of quasimaps are not easily accessible. 

Quasimap invariants of even degrees are in some sense more appealing. For example, if $\wsf=0$, then the corresponding invariants give Euler characteristics of moduli spaces of Higgs bundles. If $\wsf>0$, then by \cite[Section 7.2]{NHiggs} they determine quasimap invariants of the gerbe given by the class $\alpha$. Moreover, $K$-theoretic invariants should give more refined topological invariants. In particular, one could potentially compute these topological invariants via moduli space of degree 0 Higgs bundles  on $E$, using Vafa-Witten wall-crossing, once quasimaps to moduli space of degree 0 Higgs bundles on $E$ is better understood.
\subsection{Notation and conventions} \label{Notations} We denote the torus that scales Higgs fields by $\BC_t^*$, while the torus that scales $\p^1$ (with weight 1 at $0 \in \p^1$)  by $\BC_z^*$. We also denote 

\[ \mathbf t = \text{weight 1 representation of } \BC_t^* \text{ on } \BC,\]
\[  \mathbf z = \text{weight 1 representation of } \BC_z^* \text{ on } \BC,\]
such that $-t:=e_{\BC_t}(\mathbf t )$ and  $z:=e_{\BC_z}(\mathbf z)$ are the associated classes in the equivariant cohomology of a point. Sign in the expression for $e_{\BC_t}(\mathbf t )$ is there to match the conventions of \cite[Appendix C]{GK2}, because their torus $\BC_t^*$ scales tangent directions with weight 1, while our torus $\BC_t^*$ scales cotangent directions with weight 1. 

Moduli spaces of Higgs sheaves are not proper, hence we will always use the virtual localisation to define invariants. In order to make the notation lighter, we will denote
\[ \int_{[M]^{\mathrm{vir}}} \dots := \int_{[M^{\BC_t^*}]^{\mathrm{vir}}}\frac{\dots }{e(N^{\mathrm{vir}})}.\]

Finally, we will frequently use the fact that an obstruction theory of some space $M$ descends to the quotient $[M/G]$. That this is indeed true can be seen  either by taking quotients in the category of derived stacks, since our group actions preserve the naturally defined derived enhancements. Or we can view the descend of an obstruction theory of $M$ to $[M/G]$ as an obstruction theory of $[M/G]$ relative to $[\mathrm{pt}/G]$, which requires certain compatibility of the corresponding moduli problems, which also holds in our case. 
\subsection{Acknowledgments}  I am grateful to Sanghyeon Lee and Yaoxiong Wen for useful discussions on related topics, during which the idea of using quasisections was conceived.  I also thank Martijn Kool and Richard Thomas  for answering various questions on Vafa--Witten theory. 

This work is a part of a project that has received funding from the European Research Council (ERC) under
the European Union’s Horizon 2020 research and innovation programme (grant agreement No. 101001159).
	\section{Vafa--Witten and quasimap invariants} \label{secvwqm}
	\subsection{Preliminaries} \label{secvw}
Throughout the present work we fix a rank $\rsf$. Only in the very end of Section \ref{wallcross}, we will restrict to $\rsf=2$. We need $\rsf=2$ for the analysis of Quot schemes in Section \ref{secquot}.

	Let $C$ and $C'$ be smooth non-rational projective curves and let 
	\[L_\delta:= \CO_C(1)\boxtimes \CO_{C'}(\delta), \quad \delta \in \BQ_{>0},\]
	be an ample $\BQ$-line bundle on the product $C\times C'$. 
	 For the extremal values of $\delta$ we introduce the following notation, 
	\begin{align} \label{pm}
		\begin{split}
		\delta&=+ \quad \text{ if }\delta\gg 1, \\
		\delta&=- \quad \text{ if } \delta \ll 1. 
		\end{split}
	\end{align}
	Throughout the paper, we will be using the identification  
	\[H^2(C\times C')\cong \BZ \oplus H^1(C)\otimes H^1(C')\oplus \BZ,\]
	provided by the Kunneth's decomposition theorem. We define
	\[ \Gamma_C:=\Jac(C)[\rsf] \]
to	be a group of $\rsf$-torsion lines bundles on $C$.

		\begin{defn} \label{vwinv} Let $\dsf$ and $\asf$ be integers, such that $0 \leq  \dsf ,\asf<\rsf $. We define
	\[ M^\delta(C\times C',\dsf,\asf,\wsf) \] 
	to be a moduli space of  Higgs sheaves $(F, \phi)$ with  a fixed determinant and traceless Higgs field $\phi \in \Hom(F,F\otimes \omega_{C\times C'})$ on $C \times C'$, which are Gieseker-stable with respect to $L_\delta$. The class of $F$ is given as follows 
	\begin{align*}
		& \rk(F)=\rsf\\
	& \mathrm{c}_1(F)=(\dsf, 0, \asf)\\
	 & \Delta(F):=\mathrm{c}_1(F)^2-2\rsf\ch_2(F)=2\wsf .
	 \end{align*}
Throughout the work, we assume 
\[ \gcd(\rsf,\dsf,\asf)=1,\]
	which implies that there are no strictly semistable Higgs bundles.
	\end{defn}
\begin{rmk} The assumption on the middle component of $\mathrm{c}_1(F)$ being 0 is not restrictive, because if $C$ and $C'$ are chosen in the way that the Jacobian of one curve is not an isogenous component of another, then $H^1(C)\otimes H^1(C')$ does not contain algebraic classes. The curves $C'$ and $C$ can always be deformed to such set-up. By the deformation invariance of Vafa--Witten invariants, we can therefore assume that the middle component is zero. 
	\end{rmk}
	\subsection{Vafa--Witten wall-crossing} Conjecturally, Vafa--Witten invariants are independent of stabilities for surfaces with $p_g(S)> 1$.  Let us present more concrete evidences of this. By \cite{MM}, (physically derived) formulas for Vafa--Witten invariants for a surface $S$ with $b_1(S)=0$ and  $p_g(S)> 1$ are independent of stabilities, (see also discussion in \cite[Section 1.6]{TT}).  By \cite{DPS}, the same holds for Donaldson invariants for a surface with $b_1(S)\neq 0$ and $p_g(S)> 1$. It is therefore reasonable to expect Vafa--Witten invariants (with even insertions) for a surface with $b_1(S)\neq0$ and  $p_g(S)> 1$ are also independent of stabilities. Products of non-rational curves are among such surfaces. We now sketch a proof for this claim. 
	
	\begin{claim} \label{claim} Assume $p_g(S)>0$. If $\rsf$ is prime and there are no strictly semistable sheaves, then Vafa--Witten invariants with even $\mu$-insertions are independent of stability.  
		\end{claim}
	\textit{Sketch of Proof}. Vafa-Witten invariants consist of instanton and monopole contributions. The instanton contributions are integrals on moduli spaces of stable sheaves on the surface (descendent Donaldson invariants). On the other hand, the monopole contributions are given by integrals on moduli spaces of flags of sheaves. See \cite{TT} for more details.
	
	 The claim can therefore be proven by assembling the following results from the literature. 
\begin{itemize}
	\item Mochizuki's universal expressions for descendent invariants on moduli spaces of sheaves on surfaces, \cite{Moch}. To express instanton contributions via descendent invariants, we use G\"ottsche--Kool's expressions of the virtual equivariant Euler class in terms of descendent invariants, \cite{CY};
	\item Thomas' double-cosection argument which shows that only vertical components contribute to the monopole branch \cite[Corollary 5.30]{Th20}; 
	\item Laarakker's expressions of vertical contributions in terms of integrals on nested Hilbert schemes \cite[Theorem A]{Laa} (see also \cite[Theorem 3]{GSY} for the rank 2 case).  Laarakker's analysis extends to invariants with insertions. 
\end{itemize}
 It can be readily checked that Thomas' and Laarakker's results are independent of the assumption on $b_1(S)$. On the other hand, Mochizuki's result is more involved. 

 More conceptually, the independence of stability for Vafa--Witten invariants  should  be studied within the framework of Joyce's wall-crossing \cite{J}. 
\qed

\begin{defn}
Using the $\BC_t^*$-scaling action, we define
\begin{align*}
	&\VW^{\asf}_{\dsf,\wsf}(C\times C'):=\int_{[M^\delta(C\times C',\dsf,\asf ,\wsf)]^{\mathrm{vir}}} 1  \in \BQ 
	\end{align*}
to be Vafa--Witten invariants associated to a moduli space $ M^\delta(C\times C',\dsf,\asf,\wsf)$. By Claim \ref{claim}, they are independent of $\delta$.

\end{defn}
\subsection{Quasisection invariants} \label{squasisections} Importance of quasisections was already observed in \cite[Section 7]{Ok}. Here, we apply them in the context of a relative moduli space of Higgs bundles. 

Let $K_{C\times C'}$ be the total space of the canonical bundle $\omega_{C\times C'}$ on $C\times C'$. The variety $K_{C\times C'}$ admits projections both to $C$ and to $C'$, 
\[ \pi_C \colon K_{C\times C'} \rightarrow C, \quad   \pi_{C'} \colon K_{C\times C'} \rightarrow C'.\]
We will consider various moduli spaces (e.g.\ Quot schemes and moduli spaces of Higgs sheaves) relative to these projections. 

\begin{defn} We define $\FM^{\rel}_C(\dsf)\rightarrow C'$ to be a relative moduli space of 1-dimensional compactly supported sheaves on $\pi_{C'} \colon K_{C\times C'} \rightarrow C'$, whose associated Higgs sheaves are of rank $\rsf$, degree $\dsf$, and with a fixed determinant and a traceless Higgs field. We refer to this moduli space as a \textit{relative} moduli space of Higgs sheaves. By $M^{\rel}_C(\dsf) \rightarrow C'$ we denote its semistable locus. 
\end{defn}

Let us denote 
\[ \vsf=(\rsf,\dsf) \in H^\ev(C,\BZ). \] 
As in the absolute case, we have a determinant-line-bundle map 
\[\lambda \colon H^\ev(C,\BZ) \rightarrow \Pic(\FM^{\rel}_C(\dsf)), \]
such that a class $u \in H^\ev(C,\BZ)$ which satisfies
\[\chi(\vsf \cdot u)=\int_C\vsf \cdot u \cdot \td_C=\rsf \cdot u_2+\dsf \cdot u_1+\rsf \cdot u_1(1-g)=1\] gives a trivilisation of the  $\BC^{*}$-gerbe 
\[ \FM^\rel_C(\dsf)  \rightarrow \FM^\rel_{\mathrm{rg},C}(\dsf):=\FM^\rel_C(\dsf)\thickslash \BC^*, \]
or, in other words, a universal family on $\FM^\rel_{\mathrm{rg},C}(\dsf)$. While for  a class $u \in H^\ev(C,\BZ) $ such that $\chi(\vsf \cdot u )=0$, the line bundle $\lambda(u)$ descends to $\FM^\rel_{\mathrm{rg},C}(\dsf)$.
\\

 We define the theta line bundle $\Theta \in \Pic(\FM^{\rel}_{\mathrm{rg},C}(\dsf))$ as follows 
 \begin{align*}
 \theta&=(-\rsf,\dsf-\rsf(g-1)) \\
  \Theta&= \lambda(\theta). 
  \end{align*}
 We also define the class of $\mathrm{SL}$-trivialisations of the universal family  of $\FM^
 \rel_C(\dsf)$,
 \[ \alpha \in H^2(\FM^\rel_C(\dsf), \BZ_\rsf),\] 
 equivalently,  $\alpha$ is the Kunneth component of the first Chern class of the universal family modulo $\rsf$. The class $\alpha$ is the gerbe class of \cite{HT}.
The classes $\Theta$ and $\alpha$ will be used to define degrees of quasisections. 
\begin{rmk} If gcd$(\rsf, \dsf)$=1, then $\alpha$ is a multiple of $\Theta$ modulo $\rsf$. However, it is not the case otherwise. In any case, it is useful to keep $\alpha$ for notational purposes, because even in the case gcd$(\rsf, \dsf)$=1, invariants behave very differently depending on the degree with respect to $\alpha$.  
	\end{rmk}

\begin{defn} \label{defn} A \textit{quasisection} of $M^{\rel}_C(\dsf)$ is a section of the natural projection $p: \FM^{\rel}_{\mathrm{rg},C}(\dsf) \rightarrow C'$, 
	\[f \colon C' \rightarrow \FM^{\rel}_{\mathrm{rg},C}(\dsf), \quad p \circ f=\mathrm{id}_{C'},  \]
	 which maps generically to $M^{\rel}_C(\dsf)$.  We say that a quasisection is of degree $(\wsf,\asf) \in 
	\BZ\oplus \BZ_\rsf:=\BZ \oplus \BZ/\rsf \BZ$, if 
	 \[\deg(f^*\Theta)=\wsf, \quad f^*\alpha=\asf.\]
	 We denote the moduli space of quasisections of $M^{\rel}_C(\dsf)$ of degree $(\wsf,\asf)$ by $Q(M^{\rel}_{C}(\dsf),\asf, \wsf)$. 
	\end{defn}
\noindent The moduli spaces $Q(M^{\rel}_{C}(\dsf),\asf,\wsf)$ inherit $\BC_t^*$-actions from $\FM^{\rel}_C(\dsf)$. The properness of quasisections and existence of a perfect obstruction theory is proven in the same vain as the results of \cite{N, NHiggs}, also see \cite{LeeW} for more details. 
 \begin{defn} If $\gcd(\rsf,\dsf)=1$, we define
	\begin{align*}
		\QM^{\asf}_{\dsf, \wsf}( C)=\int_{[Q(M^{\rel}_{C}(\dsf),\asf,\wsf)]^{\mathrm{vir}}}1  \in \BQ 
	\end{align*}
	to be quasisection invariants associated to a moduli space $Q(M^{\rel}_{C}(\dsf), \asf,\wsf)$. 
\end{defn}
Notice that by a definition of a relative moduli of sheaves, a section 
\[f \colon C' \rightarrow \FM^{\rel}_C(\dsf)\] 
is given by a sheaf on 
\[K_{C\times C'}\times_{C'}C'=K_{C\times C'}. \]
Hence by \cite[Proposition 5.10]{NHiggs}, a moduli space of $L_+$-stable Higgs sheaves on $C\times C'$ is naturally a $\Gamma_{C'}$-torsor over the moduli space of quasisections  of  $M^{\rel}_{C}(\dsf)$. 
	On the other hand, the moduli space of $L_{-}$-stable Higgs sheaves on $C\times C'$  is naturally a $\Gamma_C$-torsor over the  moduli space of quasisections of $M^{\rel}_{C'}(\dsf)$. Moreover, the corresponding obstruction theories match. This is summarised in the following proposition. 
\begin{prop} \label{gmvw}  If $\gcd(\rsf,\dsf)=1$, we have
	\begin{align*}
	Q(M^{\rel}_{C'}(\dsf),\asf,\wsf)&\cong [M^{-}(C\times C',\dsf,\asf,\wsf)/\Gamma_{C}] \\ 
	Q(M^{\rel}_{C}(\dsf),\asf,\wsf)&\cong [M^+(C\times C',\dsf,\asf,\wsf)/\Gamma_{C'}],
	\end{align*}
such that the naturally defined obstruction theories on both sides match.
	\end{prop}
\textit{Proof.} Similar to \cite[Proposition 5.10]{NHiggs}, see also \cite{LeeW}. \qed
\\

\noindent We use Proposition \ref{gmvw}  as a justification for the following definition of invariants in the case of $\gcd(\rsf,\dsf)\neq 1$. 
\begin{defn} If $\gcd(\rsf,\dsf)\neq1$, we define 
	\[\QM^{\asf}_{\dsf, \wsf}( C):= \int_{[M^+(C\times C',\dsf,\asf,\wsf)/\Gamma_{C'}]^{\mathrm{vir}}}1  \in \BQ. \]
If $\gcd(\rsf,\asf)\neq1$, we define 
\[\QM^{\dsf}_{\asf, \wsf}( C):= \int_{[M^-(C\times C',\dsf,\asf,\wsf)/\Gamma_{C}]^{\mathrm{vir}}}1  \in \BQ. \]
	\end{defn}

\begin{rmk} Note that Proposition \ref{gmvw} inexplicitly depends on a choice of the universal family on the rigidified stack $\FM^{\rel}_{\mathrm{rg},C}$. We return to this point in Section \ref{hss} for an elliptic curve.
	\end{rmk}

Using Claim \ref{claim} and Proposition \ref{gmvw}, we obtain a curious correspondence between quasisection  invariants of $M^{\rel}_{C'}(\dsf)$ and $M^{\rel}_{C}(\dsf)$.
\begin{cor} \label{comp1} If $\rsf$ is prime, we have 
\begin{align*}
	&\rsf^{2g(C')}\QM^{\asf}_{\dsf, \wsf}(C)=\VW^{\asf}_{\dsf,\wsf}(C\times C')=\rsf^{2g(C)}\QM^{\dsf}_{\asf, \wsf}(C').
\end{align*}
\end{cor}
\section{Genus 1 invariants} \label{sgenus1}
\subsection{Group actions}
For the rest of the paper we assume $C'$ is an elliptic curve,
\[C'=E.\]
 Since $\pi_E \colon K_{C\times E} \rightarrow E$ is a trivial fibration, the moduli space of quasisections to $M^{\rel}_C(\dsf)$ is canonically isomorphic to a moduli space of quasimaps from $E$ to an absolute moduli space of Higgs bundles $M_C(\dsf)$ on $C$, 
 \[ Q(M^{\rel}_{C}(\dsf), \asf,\wsf) \cong Q_E(M_C(\dsf), \asf,\wsf).\]
 \noindent In fact, our primary interest is in quasimaps up to translations of $E$, i.e.\
 in the quotient
 \[[Q_{E}(M(\dsf), \asf,\wsf)/E],\] 
 where $E$ acts $Q_{E}(M(\dsf), \asf,\wsf)$ by precomposition with a translation. A similar action exists on the level of moduli spaces $M^\delta(C\times E,\dsf,\asf,\wsf)$, which we now explain. 
\\

The group 
\[E\times \Jac(E)\] 
naturally acts on sheaves.  Here, $E$ acts by pulling back a sheaf with respect to a translation $\tau_p$ by a point $p\in E$, while $\Jac(E)$ acts by tensoring a sheaf with a line bundle $L$. These operations commute. Overall, 
\[F \mapsto \tau_p^*F \otimes L.\]
Let 
\[\Phi(\CL) \subset E\times \Jac(E)\]
be the subgroup that fixes the determinant line bundle $\CL$ of sheaves in a moduli space $M^\delta(C\times E,\dsf,\asf,\wsf)$. We define 
\begin{equation} \label{phi}
 \Phi_{\asf}= (\mathrm{id},\rsf)^{-1}\Phi(\CL). 
 \end{equation}
The group $\Phi_{\asf}$ preserves rank $\rsf$ sheaves with determinant $\CL$. The action of $\Phi_{\asf}$ on sheaves therefore restricts to an action on $M^\delta(C\times E,\dsf,\asf,\wsf)$. 

By our assumption on the classes in Definition \ref{vwinv}, the line bundle $\CL$ is of the form $L \boxtimes L'$. Hence $\Phi(\CL)$ and therefore $\Phi_{\asf}$ depend only on the degree $\asf$. The group $\Phi_{\asf}$ also acts on $\FM^\rel_E(\asf)$.  In the case of 
\[ Q(M^\rel_E(\asf), \dsf,\wsf)\cong [M^-(C\times E,\dsf,\asf,\wsf)/ \Gamma_C],\]  the action of $\Phi_{\asf}$ on $Q(M^\rel_E(\asf), \dsf,\wsf)$ can be seen as identification of maps by the automorphisms of the target. 

The importance of this action is due to the next two lemmas. 
\begin{lemma} \label{canident}  There is a canonical identification
	\[ [Q_{E}(M_C(\dsf),\asf,\wsf)/E]\cong  [M^+(C\times E,\dsf,\asf,\wsf)/\Phi_{\asf}],\]
	such that the naturally defined obstruction theories on both sides match. 
	\end{lemma}
\textit{Proof.} There exists a natural map 
\begin{equation} \label{proj}
	M^+(C\times E,\dsf,\asf,\wsf) \rightarrow Q_{E}(M_C(\dsf),\asf,\wsf), 
	\end{equation}
which is a $\Gamma_E$-torsor. There also exists a natural projection
\begin{equation} \label{proj2}
\Phi_{\asf} \rightarrow E, 
\end{equation}
which is also a $\Gamma_E$-torsor. The map (\ref{proj}) is equivariant with respect to (\ref{proj2}) and the corresponding actions of $\Phi_{\asf}$ and $E$ on the source and target. It is not difficult to check, that we obtain the claimed identification after taking quotients. The rest follows from the same arguments as in \cite[Section 5.5]{NHiggs}. 
\qed
\\

\noindent The action of $\Phi_{\asf}$ can be exchanged for an insertion.  We are interested in $\mu$-insertions, which are defined as follows
\begin{align*}
 \mu \colon H^*(C\times E,\BQ) &\rightarrow H^{4-*}(H_{\BC^*}^{4-*}(\M^\delta(C\times E,\dsf, \asf, \wsf),\BQ)) \\
 \beta &\mapsto \pi_{M*}\left( \Delta(\CF)/2\rsf\cdot \pi^*_{X\times E}\beta \right ),
 \end{align*}
where $\mathsf{F}$ be the universal sheaf on $\M^\delta(C\times E,\dsf, \asf, \wsf)$. Consider now the class
\[ B_\wsf:=\frac{\mathbb 1 \boxtimes [\mathrm{pt}]}{\rsf \wsf} \in H^*(C\times E,\BQ). \]
\begin{lemma} \label{quotmu}\
	\[\VW^{\asf, \bullet}_{\dsf,\wsf}(C\times E):=\int_{[M^\delta(C\times E,\dsf,\asf,\wsf)/\Phi_{\asf}]^{\mathrm{vir}}}1=\int_{[M^\delta(C\times E,\dsf,\asf,\wsf)]^{\mathrm{vir}}}\mu(B_\wsf).  \]
	\end{lemma}
\textit{Proof.} Similar to \cite[Proposition 5.26]{NHiggs}.
\qed 
\\

By Lemma \ref{quotmu} and Claim \ref{claim}, the invariants associated to the quotient moduli space $[M^\delta(C\times E,\dsf,\asf,\wsf)/\Phi_{\asf}]$ are independent of $\delta$.

\begin{defn}
If $\gcd(\rsf,\dsf)=1$, we define 

\[\QM^{\asf, \bullet}_{\dsf, \wsf}(C)= 	\int_{[Q_E(M(\dsf),\asf,\wsf)/E]^{\mathrm{vir}}}1\in \BQ t  \]
to be invariants associated to quotient moduli spaces $[Q_E(M(\dsf), \asf,\wsf)/E]$. If $\gcd(\rsf,\asf)=1$, we also define
\[ \QM^{\dsf, \bullet}_{\asf,\wsf}( E) = \int_{[Q(M^\rel_E(\dsf), \asf,\wsf)/\Phi_{\asf}]^{\mathrm{vir}}} 1\in \BQ t\] 
to be invariants associated to the quotient moduli space $[Q(M^\rel_E(\asf), \dsf,\wsf)/\Phi_{\asf}]$. 
\end{defn} 

We use Lemma \ref{canident} as a justification for the following definition of invariants in the case of  $\gcd(\rsf,\dsf)\neq1$ and $\gcd(\rsf,\asf)$. 
\begin{defn}
	If $\gcd(\rsf,\dsf)\neq1$, we define 
\[\QM^{\asf, \bullet}_{\dsf, \wsf}(C) = \int_{[M^+(C\times E,\dsf,\asf,\wsf)/\Phi_{\asf}]^{\mathrm{vir}}}1 \in \BQ t.\]
If $\gcd(\rsf,\asf) \neq 1$, we define 
\[ \QM^{\dsf, \bullet}_{\asf, \wsf}(E)=\int_{[M^-(C\times E,\dsf,\asf,\wsf)/\Phi_{\asf}]^{\mathrm{vir}}}1 \in \BQ t.\] 
\end{defn}

\begin{rmk} Note the presence of the equivariant parameter $t$. This is due to existence of the $\BC^*_t$-equivariant cosection. We can safely remove it, thereby obtaining $\BQ$-valued invariants. 
\end{rmk}
\begin{rmk}  In the case of invariants up to translation by $E$, the role of the $E$-action is exchanged after passing from quasimaps of $M_C(\dsf)$ to quasisections of $M_E(\asf)$.  For $[Q_{E}(M(\dsf), \asf,\wsf)/E]$, taking quotient is an identification of maps by translations of the source curve $E$. On the other hand, for $[Q(M^\rel_E(\asf), \dsf,\wsf)/\Phi_{\asf}]$, taking quotient can be seen as identification of quasisections by 
automorphisms of the target $M^{\rel}_{E}(\dsf)$. 

	\end{rmk}
\noindent Using Claim \ref{claim} and Lemma \ref{quotmu}, we obtain the following result. 
\begin{cor} \label{comp} If $\rsf$ is prime, we have
	\[\QM^{\asf, \bullet}_{\dsf, \wsf}(C)=\VW^{\asf, \bullet}_{\dsf,\wsf}(C\times E)=\rsf^{2g(C)}\QM^{\dsf, \bullet}_{\asf, \wsf}(E).\] 
	\end{cor}
\begin{rmk} Note that unlike in Corollary \ref{comp1}, we don not have the factor  $\rsf^{2g(E)}$ on the left hand side. This is because the group $\Phi_\asf$ contains $\Gamma_E$. 
	\end{rmk}
\subsection{Moduli spaces of Higgs sheaves and sheaves} \label{hss} From now on, we assume that $\gcd(\rsf,\asf)=1$, unless stated otherwise.

A moduli space of rank $\rsf$ and degree $\asf$ stable Higgs $\mathrm{GL}_\rsf$-bundles on $E$, denoted by $M_E^{\mathrm{GL}}(\asf)$, is isomorphic to $K_E$ via the determinant-trace map, 
\[(\det,\tr)\colon M_E^{\mathrm{GL}}(\asf ) \xrightarrow{\cong} K_E.\] 
Hence, a moduli space of stable Higgs sheaves on $E$ with a fixed determinant and a traceless Higgs field is a point, 
\[M_E(\asf)=\{(G,0)\}=\mathrm{pt},\]
where $G$ is the unique stable sheaf with the given determinant. 
This also holds relatively for the projection $\pi_C\colon K_{C\times E} \rightarrow C$, 
\[M^{\rel}_{E}(\asf)=\{(G,0)\} \times C =C.\]
Using Corollary \ref{comp1}, we obtain an immediate consequence for quasimaps with $\wsf=0$, which confirms a part of \cite{NHiggs}[Conjecture B].
\begin{cor}  \label{const} If $\rsf$ is prime and $\asf\neq 0$, then 
	\[\QM^{\asf}_{0,0}(C)=\rsf^{2g-2}. \] 
\end{cor}

Now let us consider quasimaps with degree $\wsf \neq 0$. Since the only Higgs sheaf in $M^{\rel}_{E}(\asf)$ has a zero Higgs field, a quasisection to $M^{\rel}_{E}(\asf)$ must factor through the moduli stack of Higgs sheaves with zero Higgs fields, which is just a moduli stack of sheaves on $E$, 
\begin{equation} \label{embedding}
\FN^{\rel}_{E}(\asf) \hookrightarrow \FM^{\rel}_{E}(\asf).
\end{equation}
Since $\FN^{\rel}_{E}(\asf)$ is a relative moduli space of sheaves associated to a trivial fibration $C\times E \rightarrow C$, it  trivialises canonically 
\begin{equation} \label{trivial}
\FN^\rel_{E}(\asf)=\FN_E(\asf)\times C.
\end{equation}

\noindent Let us now consider the obstruction theory of $ \FM^{\rel}_{C}(\asf)$. Let 
\begin{equation*}\CG \in \Cohc(K_{C\times E}\times_C \FM^{\rel}_{E}(\asf)) \quad \text{and} \quad \pi\colon K_{C\times E}\times_C \FM^{\rel}_{E}(\asf) \rightarrow   \FM^{\rel}_{E}(\asf)
	\end{equation*}
 be the universal Higgs sheaf and the canonical projection. Then the virtual $C$-relative tangent complex of $\FM^{\rel}_{E}(\asf)$ is given by the following complex 
\[\BT_{\FM^{\rel}_{E}(\asf)}^{\mathrm{vir}}= R \pi_*R\CH om(\CG,\CG)_0,\] 
 where $R\CH om(\CG,\CG)_0$ is defined to be the cone
 \[ R\CH om(\CG,\CG)_0:= \mathrm{Cone}(R\pi_* R\CH om(\CG,\CG) \rightarrow H^0(K_E) \oplus H^1(\CO_E))[1]. \]
 
\noindent Note that $\BT_{\FM^{\rel}_{E}(\asf)}^{\mathrm{vir}}$ does not restrict to virtual tangent bundle of $\FN^\rel_{E}(\asf)$, hence the identification (\ref{trivial}) does not hold on the level of derived stacks.


\subsection{Chern characters} \label{chapterchern} For the purposes of wall-crossing, one needs to make a choice for a universal family on the rigidified stack $\FN_{\mathrm{rg}, E}(\asf)$. As explained in \cite[Section 3]{NHiggs}, this amounts to choosing $\usf \in H^{\mathrm{ev}}(E,\BZ)$, such that 
$\chi (\usf \cdot \vsf)=1.$

 For a choice of such $\usf=(\usf_1,\usf_2)$, the sheaf $F$ on $C\times E$ associated to a quasisection $f\colon C \rightarrow \FN^{\rel}_{E}$ of degree $\wsf$ has the following Chern character, 
\[\ch(F)=(\vsf, \check\wsf) \in  H^{\mathrm{ev}}(E,\BZ)\oplus  H^{\mathrm{ev}}(E,\BZ)(-1), \]
where $\check \wsf$ is defined by the following system of equations
	\begin{equation} \label{eqchern}
	\begin{split}
		 & \check \wsf_1\cdot \usf_2+ \check \wsf_2\cdot \usf_1=0, \\
		& \check \wsf_1  \cdot \asf -\check \wsf_2 \cdot \rsf  =\wsf.
	\end{split}
\end{equation}
For example, if $(\rsf, \asf)=(\rsf,1)$, then $\usf=(1,0)$
clearly satisfies
\[ \chi(\vsf \cdot \usf)=1. \]
Using (\ref{eqchern}), we deduce that in this case
\[ \check \wsf =(\wsf,0).\]
\section{Wall-crossing} \label{wallcross}
\subsection{$\epsilon$-stable quasisections} 
We will use Corollary \ref{comp} to compute genus 1 quasimap invariants of moduli spaces of Higgs bundles on $C$. If $\gcd(\rsf,\asf)=1$, then
\[M_E(\asf)=\{(G,0)\},\]
hence there are no sections of non-zero degree and there is a unique section of degree zero. 
 The quasimap wall-crossing for $\wsf>0$ is therefore particularly simple here, as it gives equality of invariants associated to $\epsilon=0^+$ and the wall-crossing invariants. However, there are two caveats: 
\begin{itemize}
	\item The action of $\Phi_\asf$ on $\FM^\rel_E(\asf)$,
	\item  $C$-relative set-up,
	\end{itemize}
which obscure otherwise simple computations.

Let us start with defining $\epsilon$-stable quasisections. 
From now on,  we simplify the notation in the following way, 
\begin{align*}
Q(\asf, \wsf)&:= Q(M^\rel_E(\asf), \wsf), \\
 Q(\asf, \wsf)^\bullet&:= [Q(M^\rel_E(\asf), \wsf)/\Phi_\asf],
\end{align*}
the same applies to other related spaces. 

\begin{defn} A \textit{bubbling} of $C$ is a pair $(C',\iota)$, where $C'$ a connected marked nodal curve of genus equal to $g(C)$ and
	\[\iota \colon C\hookrightarrow C'\]
	is a closed immersion.  
\end{defn}
\begin{defn}Let $\epsilon \in \BQ_{>0}$, we define $Q^{\epsilon}_{k}(\asf, \wsf)$ to be a moduli space of  quasimaps
	\[f\colon C' \rightarrow \FN^\rel_{\mathrm{rg},E}(\asf)=\FN_{\mathrm{rg},E}(\asf)\times C,\]
	such that
	\begin{itemize}
		\item $f_{ \FN(\asf)} \colon C' \rightarrow  \FN_{\mathrm{rg},E}(\asf)$ is $\epsilon$-stable (see \cite[Definition 3.5]{NHiggs}),
		\item $C'$ is a bubbling of $C$ with $k$ markings,
		\item $[f_{C}\circ \iota\colon C\rightarrow C]=\mathrm{id}_C$.
	\end{itemize}
\end{defn}

The moduli space $Q^{\epsilon}_{k}(\asf, \wsf)$ should be viewed as a moduli space of $\epsilon$-stable quasisections. The fact that these moduli spaces are proper follows from the arguments of \cite{N,NHiggs}. Recall the embedding
\[\FN^\rel_{E}(\asf) \hookrightarrow \FM^{\rel}_{E}(\asf).\] 
With respect to this embedding, we endow $Q^{\epsilon}_{n}(\asf, \wsf)$ with the obstruction theory given by the complex
\[\mathbb T^{\mathrm{vir}}_{Q^{\epsilon}}:=\pi_*f^*\BT_{\FM^{\rel}_{E}}^{\mathrm{vir}}. \]
Its perfectness is proven in the same vain as in \cite[Proposition 3.12]{NHiggs}.

Let us indicate what moduli spaces $Q^{\epsilon}_{k}(\asf, \wsf)$ are for the extremal values of $\epsilon$. Using the same notation as in (\ref{pm}), if $\epsilon=-$ , we get 
\[Q^{-}_{0}(\asf, \wsf) \cong Q(\asf, \wsf).  \]
If $\epsilon=+$, then
\begin{align*}Q^{+}_{0}(\asf, \wsf)&=\mathrm{pt} \quad \text{if } \wsf=0, \\
	Q^{+}_{0}(\asf, \wsf)&=\emptyset \quad \hspace{0.17cm} \text{if } \wsf\neq0.
\end{align*}

\noindent We now discuss the wall-crossing between invariants associated to different values of $\epsilon$. 
\begin{defn}
	Let $GQ(\asf, \wsf)$ be a moduli space  of prestable quasimaps   
	\[ f\colon \p^1 \rightarrow \FN^\rel_{E}(\asf)=\FN_{E}(\asf)\times C. \]
	The $\BC^*_z$-action on the source $\p^1$ induces a $\BC^*_z$-action  on $GQ(\asf, \wsf)$. We define 
	\[W^\rel(\asf,\wsf) \subset  GQ(\asf, \wsf)\]
	to be the component of the $\BC^*_z$-fixed locus of quasimaps with base points only at $0\in \p^1$.  
\end{defn}
We endow  $GQ(\asf, \wsf)$  with obstruction theory given by the complex 
\[\BT^{\mathrm{vir}}_{GQ}:=\pi_*f^*\BT_{\FM^{\rel}_{E}}^{\mathrm{vir}},\] 
using the embedding $\FN^\rel_{E}(\asf) \hookrightarrow \FM^{\rel}_{E}(\asf)$. The space $ W^\rel(\asf,\wsf)$ inherits  the obstruction theory defined by the fixed part of the obstruction theory of $GQ(M^\rel_E(\asf), \wsf)$, as well as the virtual normal bundle $N^{\mathrm{vir}}$ defined by the moving part of the obstruction theory. 
\subsection{Moduli spaces of flags}
As before, there exists a canonical identification of moduli spaces
\begin{equation} \label{ident}
	W^\rel(\asf,\wsf)=W(\asf,\wsf) \times C
\end{equation}
where $W(\asf,\wsf)$ is analogous space defined via quasimaps to $\FN_E(\asf)$.  By \cite[Section 4]{Ob}, the space $W(\asf,\wsf)$ admits a description in terms of moduli spaces of flags, which we now recall. In what follows by $G$ we denote the sheaf supported on the zero section in $K_E$. Also, we define 
\[ \mathbf z := \text{weight 1 representation of } \BC_z^* \text{ on } \BC.\]

\noindent Let
\[\mathrm{Fl}(\asf)=\{ F_1 \subseteq F_2 \subseteq \dots \subset F_{r-1} \subseteq F_r=G \}\]
be a moduli space of  flags, such that consecutive terms are allowed to be equal. To each  $F_{\bullet}$ and a choice of an integer $k\in \BZ$, we can associated a $\BC^*$-equivariant  torsion-free sheaf $\CF$ on $K_E\times \BA^1$, 
\begin{equation}
	\CF= F_1 \mathbf z^{k+1} \oplus F_2 \mathbf z^{k+2}  \oplus \dots F_{r-1} \mathbf z^{k+r-1} \oplus G \mathbf z^{k+r} \oplus G \mathbf z^{k+r+1} \dots 
\end{equation}
In fact,  for torsion-free sheaves $G$, such association is an equivalence between $\BC^*$-fixed torsion-free sheaves on $K_E\times \BA^1$ and weighted flags (up to a choice of $k$).  Moreover, 
each $\BC^*$-fixed sheaf in $W(\asf,\wsf)$ is canonically $\BC^*$-equivariant. Let us denote  by
\[ \mathrm{Fl}(\asf,\wsf) \subset \mathrm{Fl}(\asf)\] 
the locus of flags that correspond to sheaves in $W(\asf,\wsf)$. By construction, we have 
\begin{equation} \label{flagsf}
	W(\asf,\wsf) \cong \mathrm{Fl}(\asf,\wsf). 
\end{equation}
Analogously,  let  $\mathrm{Fl}^\rel(\asf,\wsf)$ be the relative moduli space of flags of $(\pi^*_EG)_{|C\times E}$ on the relative surface $\pi_C \colon K_{C\times E} \rightarrow C$. Quotients of $G$ viewed as quotients of a sheaf on $E$, on $K_E$ or on  $K_{C\times E}$ are the same. Hence identification (\ref{flagsf}) also holds relatively, 
\begin{equation} \label{identrel}
	W^\rel(\asf,\wsf)  \cong \mathrm{Fl}^\rel(\asf,\wsf)\cong  \mathrm{Fl}(\asf,\wsf) \times C.
\end{equation}

\subsection{Obstruction theory of flags} By \cite[Section 4]{Ob}, obstruction theories of moduli spaces $W^\rel(\asf,\wsf)$ and $ \mathrm{Fl}^\rel(\asf,\wsf)$ also agree. We now describe the obstruction theory of $W^\rel(\asf,\wsf)$ and the associated virtual normal bundle $N^{\mathrm{vir}}$.  Let 
\[\CF_1 \subseteq \CF_2\subseteq \dots \subseteq \CF_r=G \]
be the universal flag on $K_{C\times E} \times_C \mathrm{Fl}^\rel(\asf,\wsf)$ and
\[  \pi\colon K_{C\times E} \times_C \mathrm{Fl}^\rel(\asf,\wsf) \rightarrow \mathrm{Fl}^\rel(\asf,\wsf) \] 
be the natural projection.
\begin{thm} \label{obstth}
	
The obstruction theory of $\mathrm{Fl}^\rel(\asf,\wsf)$ is given by the complex
\[ \mathbb T^{\mathrm{vir}}_{\mathrm{Fl}} =\mathrm{Cone} \left( \bigoplus^{i=r-1}_{i=1} R\CH om_\pi(\CF_i, \CF_i) \rightarrow \bigoplus^{i=r-1}_{i=1}R\CH om_\pi(\CF_i, \CF_{i+1}) \right).\]
The $K$-class of $N^{\mathrm{vir}}$ is 
\begin{align*} \label{Euler}
	N^{\mathrm{vir}}= &-\sum_{i\geq1}\sum_{k\geq1} R \CH om_\pi(\CF_{i+k}/\CF_{i+k-1},\CF_i) \mathbf z^{-k} \\
	&+ \sum_{i\geq1}\sum_{k\geq1}   R \CH om_\pi(\CF_{i+k+1}/\CF_{i+k},\CF_i)^\vee\mathbf z^{k}.
	\end{align*}
\end{thm}
\textit{Proof.} See \cite[Section 4]{Ob}. 
\qed 
\\

We will need the following technical lemma for later. 
\begin{lemma} 	\label{lemma}  We have the following identity in the K-group 
	\[R \CH om_\pi(\CF_{j+1}/\CF_{j},\CF_i)=\mathsf{K}(1- \omega_C\mathbf t) \]
	for some $K$-class $\mathsf{K}$.
	\end{lemma}

\textit{Proof.} Let us denote $\CA:=\CF_{j+1}/\CF_j$ and $\CB:=\CF_{i}$. Both $\CA$ and $\CB$ are scheme-theoretically supported on the zero section $C\times E  \subset K_{C\times E}$, they can therefore be extended to the entire $K_{C\times E}$ by pulling them back by the projection $K_{C\times E} \rightarrow C\times E$. We denote these extensions by $\bar \CA$ and $\bar{\CB}$. 

Consider now the sequence on $K_{C\times E}$,
\[ 0 \rightarrow \CO(-C\times E) \rightarrow  \CO_{K_{C\times E}} \rightarrow \CO_{C\times E} \rightarrow 0,\]
we  tensor it with $\bar \CA$, 
\[0 \rightarrow \bar \CA(-C\times E) \rightarrow \bar \CA \rightarrow \CA \rightarrow 0,\]
then  we apply $ R \CH om_\pi( -, \CB)$ to obtain  the distinguished triangle 
\begin{equation} \label{triangle}
  R \CH om_\pi( \CA,  \CB) \rightarrow   R \CH om_\pi( \bar \CA, \CB) \rightarrow R \CH om_\pi( \bar \CA(-C\times E), \CB) \rightarrow. 
 \end{equation}
There is a natural $\BC^*_t$-equivariant identification 
\[\CO_{K_{C\times E}}(-C\times E)_{|C\times E} \cong \omega^{\vee}_{C\times E} \mathbf{t}^{-1} \cong \omega^{\vee}_C \mathbf t^{-1},\]
which gives us that
\[ R \CH om_\pi( \bar \CA(-E), \CB) \cong R \CH om_\pi( \bar \CA, \CB) \boxtimes \omega_C\mathbf t.\]
Passing to the $K$-group, the distinguished triangle (\ref{triangle}) therefore gives us that
\[ R \CH om_\pi( \CA,  \CB)=R \CH om_\pi( \bar \CA, \CB)(1- \omega_C\mathbf t).\]
\qed

\subsection{Cosections} \label{scosections}As always, the identification 
\[\mathrm{Fl}^\rel(\asf,\wsf)\cong  \mathrm{Fl}(\asf,\wsf) \times C\] 
does not respect the obstruction theory, as there is a twist coming from $C$, which becomes visible on the level of cosections. We now remind the reader of these cosections. Let 
\[ \mathrm{Quot}^\rel(\asf,\wsf) \subset \mathrm{Fl}^\rel(\asf,\wsf)\] 
be the locus of Quot schemes, i.e.\ one-step flags. By 
\[ \mathrm{Quot}^\rel(\asf,\wsf)^{c} \subset \mathrm{Fl}^\rel(\asf,\wsf) \]
we denote its complement.  Cosections are constructed in exactly the same way as in \cite[Section 5.4]{KKV} (see also \cite[Section 4]{Ob} and \cite[Section 11.2]{N}). The relative canonical sheaf of 
\[K_{C\times E} \rightarrow C\]
is the pullback of $\omega_C^\vee$. The argument of \cite[Section 11.2]{N}, which uses Serre's duality, shows that we now have a cosection to $\omega_C^{ \oplus 2}$ (here, we use the identification (\ref{ident})) instead of the trivial line bundle $\CO^{\oplus 2}$, 

\begin{equation*}
	\sigma=(\sigma_1, \sigma_2) \colon h^1(\BT_{\mathrm{Fl}}^{\mathrm{vir}}) \rightarrow  \omega_C^{ \oplus 2}\mathbf t,
\end{equation*}
where 
\[ \mathbf t := \text{weight 1 representation of } \BC_t^* \text{ on } \BC.\]
As in \cite[Proposition 10.6]{N}, we have the following result. 
\begin{prop} \label{cosections}
	The twisted cosection $\sigma$ is surjective on $\mathrm{Quot}^\rel(\asf,\wsf)^{c} $ in  $W^\rel(\asf,\wsf)$. On $\mathrm{Quot}^\rel(\asf,\wsf)$, only the component $\sigma_1$ is surjective. 
\end{prop}
\textit{Proof} Similar to \cite[Proposition 12]{KKV}. \qed 
\\

By the description of the obstruction theory of $W^\rel(\asf,\wsf)$ in terms of flags from \cite[Section 4]{Ob}, we can compute its virtual dimension. Indeed, for any two sheaves $F_1$ and $F_2$ supported on the zero section of $K_E$, we have 
\begin{equation} \label{dimension}
\sum \mathrm{ext}^i(F_1 ,F_2) =\ch(F_1^\vee) \cdot \ch(F_2) =0,
\end{equation}
the $C$-relative virtual dimension of $\mathrm{Fl}^\rel(\asf,\wsf)$ is therefore 0. Hence the absolute virtual dimension of $\mathrm{Fl}^\rel(\asf,\wsf)$ is 1. 

Both the virtual normal bundle and the cosections are $\Phi_\asf$-equivariant by the construction, hence they descend to the quotients 
\[[\mathrm{Fl}^\rel(\asf,\wsf)/\Phi_\asf].\]
This, in conjunction with Proposition \ref{cosections}, implies that the virtual fundamental cycles of $[\mathrm{Fl}^\rel(\asf,\wsf)/\Phi_\asf]$, whence restricted to Quot schemes and to their complements, are of the following form.
\begin{cor}\label{classes}
\begin{align}  \
	\begin{split}
		[\mathrm{Quot}^\rel(\asf,\wsf)/\Phi_\asf]^{\mathrm{vir}}&=\mathsf{B}\boxtimes[\mathrm{pt}] \in H^{}_0(W(\asf,\wsf) \times C,\BQ)[t] \\
		[\mathrm{Quot}^\rel(\asf,\wsf)^c/\Phi_\asf]^{\mathrm{vir}}&=t\mathsf{B}'\boxtimes[\mathrm{pt}] \in H_2(W(\asf,\wsf) \times C,\BQ)[t].
	\end{split}
\end{align}

\end{cor}

\subsection{Master space}
Let $\epsilon^0 \in \BQ_{>0}$ be a wall of $\epsilon$-stabilities for quasisections, let $\epsilon^+$ and $\epsilon^-$ be the values close to the wall $\epsilon^0$ from the right and from the left respectively. Consider the master space $MQ^{\epsilon_0}(\asf, \wsf)$ for the wall-crossing around the wall $\epsilon_0$ (we refer to \cite[Section 4]{YZ} for the construction of the master space). Let $\wsf_0=1/\epsilon_0$.  By construction, there is a $\BC_z^*$-action on $MQ^{\epsilon_0}(\asf, \wsf)$. In what follows we use the identification
\[M^{\rel}_{E}(\asf)=\{(G,0)\} \times C =C.\]

 Let us define the following relation on spaces.
\[ X \approx Y \iff X=Y \textit{ up to \'etale covers, finite gerbes and virtual blow-ups.}\]
Any of the operations above have only a mild affect on the virtual invariants and for the purpose of our arguments they can be ignored, see \cite{YZ} for more details. 
The $\BC_z^*$-fixed locus of $MQ^{\epsilon_0}(\asf,\wsf)$ then has the following form, 
\begin{equation*}
MQ^{\epsilon_0}(\asf, \wsf)^{\BC_z^*}  
\approx Q^{\epsilon^+}(\asf, \wsf) \cup Q^{\epsilon^-}(\asf, \wsf)\cup \coprod_{k } \left(Q_{k}^{\epsilon^+}(\asf, \wsf') \times_{C^k}  W^\rel(\asf,\wsf_0)^k \right),
\end{equation*}
such that $\wsf=\wsf'+k  \wsf_0$.

The group $\Phi_\asf$ acts on the master space $MQ^{\epsilon_0}(\asf, \wsf)$.  Since the action of $\Phi_\asf$ and $\BC_z^*$ on $MQ^{\epsilon_0}(\asf, \wsf)^{\BC_z^*}$ commute, operations of taking quotient by $\Phi_\asf$ and taking $\BC_z^*$-fixed locus also commute, we therefore obtain 
\begin{multline} \label{fixedlocus}
	[MQ^{\epsilon_0}(\asf, \wsf)^{\BC_z^*}/\Phi_\asf] \\
	\approx Q^{\epsilon^+}(\asf,\wsf)^\bullet \cup Q^{\epsilon^-}(\asf, \wsf)^\bullet\cup \coprod_{k} \left( \left[ Q_k^{\epsilon^+}(\asf, \wsf') \times_{C^k}  W^\rel(\asf,\wsf_0)^k/\Phi_\asf\right] \right),
\end{multline}
such that the action on the wall-crossing components (components on the right in the expression above) is given by the diagonal action of $\Phi_\asf$. By \cite[Section 6]{YZ}, the wall-crossing formula is obtained by taking residues of the localisation formula associated to (\ref{fixedlocus}).  Let $\CN^{\mathrm{vir}}$ be the virtual normal bundle of wall-crossing components.  The wall-crossing invariants are therefore given by the following residues (up to some factors coming from the relation "$\approx$").

\[ \mathrm{Res}_{z=0}\ \left( \frac{\left[ Q_k^{\epsilon^+}(\asf, \wsf') \times_{C^k}  W^\rel(\asf,\wsf_0)^k/\Phi_\asf\right]}{e_{\BC_{z,t}^*}(\CN^{\mathrm{vir}})} \right).\]
We now show that most of the wall-crossing invariants vanish essentially by the second cosection arguments, except that our cosections are twisted as explained in Section \ref{scosections}, which forces us to work a bit harder to obtain vanishing. 
\begin{thm} \label{thetheorem} If $\epsilon_0=1/\wsf$, then
	\begin{multline*}
		\deg[Q^{\epsilon^+}(\asf, \wsf)^\bullet]^{\mathrm{vir}} - \deg [Q^{\epsilon^-}(\asf, \wsf)^\bullet]^{\mathrm{vir}} = \deg \mathrm{Res}_{z=0} \left(  \frac{[\mathrm{Quot}^\rel(\asf,\wsf)/\Phi_\asf]^{\mathrm{vir}}}{e_{\BC_{z,t}^*}(N^{\mathrm{vir}})}  \right)
		\end{multline*}
	Otherwise, 
	\[\deg[Q^{\epsilon^+}(\asf, \wsf)^\bullet]^{\mathrm{vir}} =\deg [Q^{ \epsilon^-}(\asf, \wsf)^\bullet]^{\mathrm{vir}}.\]
\end{thm} 



\subsection{Proof of Theorem \ref{thetheorem}}
By \cite[Section 6]{YZ}, we just have to analyse the wall-crossing components, which we now will do (see also \cite[Section 6]{N}, \cite[Section 10]{NHiggs} and \cite{LeeW}). 
\subsubsection{}
Assume $k \geq 2$ or $\wsf'\neq 0$, then the space 
\[X:=[Q_k^{\epsilon^+}(\asf, \wsf') \times_{C^k} W^\rel(\asf,\wsf_0)^k/\Phi_\asf]\]
have an additional action of $\Phi_\asf$ coming from any of the components of the product. To distinguish it from the diagonal action of $\Phi_\asf$, we denote it by $\Phi_\asf'$. The obstruction theory of the quotient  $[X/\Phi_\asf']$ is compatible with the obstruction theory of $X$, hence 
\[ \pi^*[X/\Phi_\asf']^{\mathrm{vir}}=[X]^{\mathrm{vir}}\]
for the quotient map $\pi \colon  X\rightarrow  [X/\Phi_\asf'] $.  Moreover, the virtual normal bundles of \cite{YZ} are $\Phi_\asf'$-equivariant, hence descend to  the quotient $[X/\Phi_\asf']$. Overall, we obtain that the wall-crossing class is a pullback of some class $\mathsf{A}$ from the quotient $[X/\Phi_\asf']$, 
\[ \frac{\left[ Q_k^{\epsilon^+}(\asf, \wsf') \times_{C^k}  W^\rel(\asf,\wsf_0)^k/\Phi_\asf\right]}{e_{\BC_{z,t}^*}(\CN^{\mathrm{vir}})}=\pi^*\mathsf{A},\]
its degree is therefore 0 and it does not contribute to the wall-crossing formula,
\[\deg \mathrm{Res}_{z=0} \left( {\frac{\left[ Q_k^{\epsilon^+}(\asf, \wsf') \times_{C^k} W^\rel(\asf,\wsf_0)^k/\Phi_\asf\right]}{e_{\BC_{z,t}^*}(\CN^{\mathrm{vir}})}} \right)=0. \] 
\subsubsection{} It remains to determine the contribution of terms 
\[ [Q_{1}^{ +}(\asf, 0)\times_C  W^\rel(\asf,\wsf)/\Phi_\asf]. \]
Firstly, $Q_{1}^{ +}(\asf, 0)=C$, hence 
\[ [Q_{1}^{ +}(\asf, 0)\times_C  W^\rel(\asf,\wsf)/\Phi_\asf]=[W^\rel(\asf,\wsf)/\Phi_\asf].\]
We will now show that the complement of $\mathrm{Quot}^\rel(\asf,\wsf) \subset W^\rel(\asf,\wsf)$ does not contribute. Ideally, one would say that this follows from the double cosection argument. However, in this case, the cosections are twisted due to the relative set-up, hence one has to do a little bit of work. By \cite{YZ}, 
\begin{equation} \label{eq}
	\mathrm{Res}_{z=0}\left(\frac{[W^\rel(\asf,\wsf)/\Phi_\asf]^{\mathrm{vir}}}{e_{\BC_{z,t}^*}(\CN^{\mathrm{vir}})}\right)
	= \mathrm{Res}_{z=0} \left( \frac{[W^\rel(\asf,\wsf)/\Phi_\asf]^{\mathrm{vir}}}{e_{\BC_{z,t}^*}(N^{\mathrm{vir}})} \right).
\end{equation}
Note that higher powers of $z$ do not contribute by the dimension constraint,  as $\dim Q_{1}^{ +}(\asf, 0)=1$.

We argue that $ \mathrm{Quot}^\rel(\asf,\wsf)^c \subset W^\rel(\asf,\wsf)$ does not contribute, because the quantity (\ref{eq}) is a multiple of $t^2$. As taking quotient by $\Phi_\asf$ can be exchanged with taking an insertion, it is enough to show that (\ref{eq}) is a multiple of $t^2$ before taking quotient.  By Corollary \ref{classes}, we know that $ [ \mathrm{Quot}^\rel(\asf,\wsf)^c]^{\mathrm{vir}}$ is a multiple of $t$, we therefore have to show that residue of $e_{\BC_{z,t}^*}(N^{\mathrm{vir}})^{-1}$ is a multiple of $t$ too. By Theorem \ref{obstth}, the class $e_{\BC_{z,t}^*}(N^{\mathrm{vir}})^{-1}$ admits the following expression
\begin{align} \label{Eulerclass}
 e_{\BC_{z,t}^*}(N^{\mathrm{vir}})^{-1}&= \prod_{i\geq 1, k\geq 1} 
\frac{e_{\BC_{z,t}^*} ( R \CH om(\CF_{i+k}/\CF_{i+k-1},\CF_i) \mathbf z^{-k})}{e_{\BC_{z,t}^*} ( R \CH om(\CF_{i+k+1}/ \CF_{i+k},\CF_i)^\vee  \mathbf z^{k})} .
\end{align}
By (\ref{dimension}), the rank of $R \CH om(\CF_{i+k}/ \CF_{i+k-1},\CF_i)$  is 0, hence
\begin{multline*}e_{\BC_{z,t}^*} ( R \CH om(\CF_{i+k} /\CF_{i+k-1},\CF_i) \otimes \mathbf z^{-k}) \\
	=\sum_{j\geq0}(kz)^{-j}\mathrm{c}_j(R \CH om(\CF_{i+k} /\CF_{i+k-1},\CF_i)) \\	
	=1-(kz)^{-1}\mathrm{c}_1(R \CH om(\CF_{i+k+1} /\CF_{i+k},\CF_i))+ O(z^{-2}),
	\end{multline*}
the same applies to the denominator of (\ref{Eulerclass}). 
We therefore obtain that 
\begin{multline} \label{O}
 e_{\BC_{z,t}^*}(N^{\mathrm{vir}})^{-1} = 1-\sum_{i,k} (kz)^{-1}\mathrm{c}_1(R \CH om(\CF_{i+k} /\CF_{i+k-1},F_i))\\
 -  \sum_{i,k} (kz)^{-1}\mathrm{c}_1(R \CH om(\CF_{i+k+1} /\CF_{i+k},\CF_i))^\vee + O(z^{-2}).
 \end{multline}
By Lemma \ref{lemma}, we obtain that 
\[\mathrm{c}_{1} ( R \CH om(\CF_{j}/ \CF_{j-1}, \CF_i)=  \mathsf{A} t + \mathsf{A}\cdot \mathrm{c}_1( \omega_C ) \in H^2(\mathrm{Quot}^c,\BZ),\]
for some class $\mathsf{A}$ of cohomological degree 0. 
Using Corollary \ref{classes} and (\ref{O}), we conclude that
\[\mathrm{Res}_{z=0} \left( \frac{[\mathrm{Quot}^\rel(\asf,\wsf)^c]^{\mathrm{vir}}}{e_{\BC_{z,t}^*}(N^{\mathrm{vir}})} \right)= \mathsf{A}' t^2\in H_2(\mathrm{Quot}^c,\BZ)[t], \]
for some class $\mathsf{A}'$ of homological degree 2. Taking degree of the class above, we therefore obtain 0. The finishes the proof of Theorem \ref{thetheorem}. 

\subsection{Contributions from Quot schemes} \label{contributions}

We now have to determine the contributions of $\mathrm{Quot}^\rel(\asf,\wsf) \subset W^\rel(\asf,\wsf)$.  Firstly, by \cite[Section 4]{Ob}, the component $\mathrm{Quot}^\rel(\asf,\wsf) \subset W^\rel(\asf,\wsf)$ is composed of the following Quot schemes, 
\[\mathrm{Quot}^\rel(\asf,\wsf) = \coprod_{\msf | \wsf}\mathrm{Quot}^\rel(\asf, u_\msf),\]
the classes $u_\msf$ are  defined as follows
\[u_\msf:=h_\msf\vsf-\frac{\check{\wsf}}{\msf},\]
where $\check \wsf$ is given by $(\ref{eqchern})$ and $h_\msf$ is the unique integer such that 
\[h_\msf\rsf -\frac{\check{\wsf}_1}{\msf}\in [0, \rsf-1].\]
 We therefore obtain the following proposition.
\begin{prop}\label{wallclass}  \
\[	\deg [\mathrm{Quot}^\rel(\asf,\wsf)/\Phi_\asf]^{\mathrm{vir}}= \sum_{\msf | \wsf} \deg [\mathrm{Quot}^\rel(\asf,u_\msf)/\Phi_\asf]^{\mathrm{vir}}.\]
\end{prop}
\textit{Proof.} See \cite{Ob}[Section 4]. \qed
\\

Let us analyse $N^{\mathrm{vir}}$ over each component $[\mathrm{Quot}^\rel(\asf,u_\msf)/\Phi_\asf ]$. In what follows we use the notations from Section \ref{secquot}. By \cite[Section 4.4]{Ob}, the equivariant Euler class of the virtual normal bundle $N^{\mathrm{vir}}$ can be expressed as follows
\begin{align*}
	e_{\BC_{z,t}^*}(N^{\mathrm{vir}})^{-1}&=e_{\BC_{z,t}^*}(R \CH om_\pi (\CQ^\rel, \CK^\rel) \mathbf z^{-\msf}) \\
	&=e_{\BC_{z,t}^*}(R \CH om_
	\pi(\CK^\rel,\CQ^\rel)^{\vee}\mathbf t \mathbf  z^{-\msf})\\
	&=\sum_{k \in \BZ}(-\msf z)^{-k}\mathrm{c}_{k}(R\CH om_\pi(\CK^\rel,\CQ^\rel)^{\vee}\mathbf t),
\end{align*}
where $\CK^\rel$ and $\CQ^\rel$ are as in Section \ref{secquot}. We are interested in the residue of  $e_{\BC_{z,t}^*}(N^{\mathrm{vir}})^{-1}$, 
\begin{align*}
	\mathrm{Res}_{z=0}(e_{\BC_{z,t}^*}(N^{\mathrm{vir}})^{-1})=& -\msf^{-1}\mathrm{c}_{1}(R\CH om_\pi(\CK^\rel,\CQ^\rel)^{\vee}\mathbf t) \\
	=&-\msf^{-1} \rk(\CH om_\pi(\CK^\rel,\CQ^\rel)^{\vee})(\mathrm{c}_1(\omega_C)-t) \\
	=&-\msf^{-1}\dim(\mathrm{Quot}(\asf,u_\msf))(\mathrm{c}_1(\omega_C)-t),
\end{align*}
where the second equality follows from Lemma \ref{relquot} and (\ref{tangent}). Also recall from Section \ref{Notations} that by convention we set $-t:=e_{\BC_t}(\mathbf t )$. Using Corollary \ref{classes}, the total residue then takes the following form
\begin{multline*}
\mathrm{Res}_{z=0}\left(\frac{[\mathrm{Quot}^\rel(\asf, u_\msf )/\Phi_\asf]^{\mathrm{vir}}}{e_{\BC_{z,t}^*}(N^{\mathrm{vir}})}\right) \\
=\msf^{-1}\dim(\mathrm{Quot}(\asf,u_\msf))  [ \mathrm{Quot}^\rel(\asf, u_\msf )/\Phi_\asf]^{\mathrm{vir}}t.
\end{multline*}

Assume that $(\rsf,\asf)=(2,1)$, using the analysis from Corollary \ref{lemma1} and Corollary \ref{lemma2}, we get 
\[ \deg \mathrm{Res}_{z=0}\left(\frac{[\mathrm{Quot}^\rel(\asf, u_\msf )/\Phi_\asf]^{\mathrm{vir}}}{e_{\BC^*}(N^{\mathrm{vir}})}\right)=(2g-2)\msf^{-1}t .\]
Now, applying Theorem \ref{thetheorem} repeatedly and using the fact that $Q^{+}(\asf, \wsf)$ is empty for $\wsf\neq 0$, we obtain the following result (for how $\dsf$ and $\wsf$ are related, see (\ref{eqchern})).
\begin{thm} \label{beforemain} If $(\rsf,\asf)=(2,1)$,  then
\label{wsE2}
\begin{align*} 
	\QM(E)_{1, \wsf}^{\dsf, \bullet}=\begin{cases}
		(2g-2) \sum_{\msf| \wsf}\msf^{-1}t, \quad &\text{ if } \wsf=\dsf \ \mathrm{mod}\ 2 \\ 
		0, &\text{ otherwise.}
		\end{cases}
\end{align*}
\end{thm}
\noindent Using Corollary \ref{comp}, we obtain the desired quasimap invariants. 
\begin{thm} \label{themain} If $(\rsf,\asf)=(2,1)$, then
	\[ \QM(C)_{\dsf,\wsf}^{1, \bullet} = \begin{cases}
		(2g-2)2^{2g} \sum_{\msf| \wsf}\msf^{-1}t, \quad &\text{ if } \wsf=\dsf \ \mathrm{mod}\ 2 \\ 
		0, &\text{ otherwise.}
		\end{cases} \]
	\end{thm}
\noindent This gives us  Theorem \ref{theeq2} and \ref{theeq} (after passing to reduced invariants, i.e.\ after dividing by $t$).
\begin{rmk} Since $M_E(1)=\{(G,0)\}$,  invariants $\QM(C)_{\dsf,\wsf}^{1, \bullet} $ have only instanton contributions, i.e.\  on $C$ they correspond to invariants of  moduli spaces of stables bundles $T^*N_C(\dsf) \subset M_C(\dsf)$. On the contrary, the even-degree invariants are completely monopole, i.e.\ they correspond to invariants of the complement of $N_C(\dsf)$ in the nilpotent cone. This was expected from \cite{MM}, see also \cite[Remark 7.2]{NHiggs}. 
\end{rmk}
 \begin{rmk}  \label{genus1}Let us now comment on the fact that genus 1 Gromov-Witten invariants of $E$ have very similar expressions, as it was mentioned in introduction. If $\gcd(\rsf,\dsf)=1$, then a moduli space of stable sheaves on $E$ is naturally isomorphic to $E$ via the determinant map. Hence Gromov--Witten theory of $E$ is equivalent to the one of its moduli spaces of sheaves. Here, we study some kind of twisted Gromov--Witten theory of moduli spaces of sheaves on $E$.  Hence (posteriori) it is not so unexpected that we get similar answers. Perhaps this phenomenon can even be made precise.  
\end{rmk}
\subsection{Higher rank}
By (\ref{eqchern}), if we assume that all divisors $\msf$ of $\wsf$ are congruent to $0$ or $\asf$ modulo $\rsf$, then 
\[u_\msf=(0,\ksf) \ \text{or} \ (\rsf-1,\ksf),\]
therefore using the arguments of the previous section, the analysis of Section \ref{relevant} is enough to conclude that following. 

\begin{prop} \label{themain2} If $\rsf$ is prime,  $\asf\neq 0$ and all divisors $\msf$ of $\wsf$ satisfy 
	\[  \msf= 0 \ \mathrm{ or } \ \asf \ \mathrm{ mod } \ \rsf,\]
	   then
		\[ \QM(C)_{\dsf,\wsf}^{\asf, \bullet} = \begin{cases}
			(2g-2)\rsf^{2g} \sum_{\msf| \wsf}\msf^{-1}t, \quad &\text{ if } \wsf= \dsf \cdot \asf \ \mathrm{mod}\ \rsf \\ 
			0, &\text{ otherwise.}
		\end{cases} \]
	
	\end{prop}

This agrees with \cite[Conjecture E]{NHiggs}. Our methods involving Quot schemes lead to an obvious conjectural extension of Proposition \ref{themain2}.
\begin{conjecture} \label{conj}If $\gcd(\rsf,\asf)=1$, then
	\begin{align*} 
		\QM(C)_{\asf, \wsf}^{\dsf, \bullet}=\begin{cases}
			(2g-2) \rsf^{2g}\sum_{\msf| \wsf}\msf^{-1}t, \quad &\text{ if } \wsf=\dsf \cdot \asf \ \mathrm{mod}\ \rsf \\ 
			0, &\text{ otherwise.}
		\end{cases}	
	\end{align*}
	
	\end{conjecture} 

\section{Quot schemes} \label{secquot}

\subsection{Relative Quot schemes} 

Let $\mathrm{Quot}^\rel(\asf)$ be the relative Quot associated to the relative surface
\[\pi_C \colon K_{C\times E} \rightarrow C. \]
and the sheaf $(\pi^*_EG)_{|C\times E}$ supported on the zero section $C\times E  \hookrightarrow K_{C\times E}$. Here, the sheaf $G$ is the unique stable sheaf with a fixed determinant of rank $\rsf$ and degree $\asf$ on $E$. 
Considering quotients of $(\pi^*_EG)_{|C\times E}$ as sheaves  on $K_{C\times E}$ or as sheaves on $C\times E$ is equivalent. Moreover, since $C\times E \rightarrow C$ is a trivial fibration and $(\pi^*_EG)_{|C\times E}$ is scheme-theoretically supported on $C\times E$,  we have a canonical identification
\begin{equation} \label{ident1}
	\mathrm{Quot}^{\rel}(\asf)=\mathrm{Quot}(\asf) \times C,
\end{equation}
where $\mathrm{Quot}(\asf)$ is the Quot scheme of $G$ on $E$. However, the natural obstruction theories are not respected by the identification (\ref{ident1}), even if we endow $\mathrm{Quot}(\asf)$ with the obstruction theory associated to 1-dimensional sheaves on $K_E$, because there is a twist coming from $C$.  Let us now elaborate on this point by comparing the obstruction theories in both cases. 
\subsection{Relative vs.\ absolute}
Let $\CK$ and $\CQ$ be the universal kernel and quotient sheaves on $K_E \times \mathrm{Quot}(\asf)$, and let $\pi \colon K_E \times \mathrm{Quot}(\asf) \rightarrow \mathrm{Quot}(\asf)$ be the natural projection. The $\BC_t^*$-equivariant obstruction theory of $\mathrm{Quot}(\asf)$ is given by the complex
\begin{equation} \label{complex}
	\BT_{\mathrm{Quot}}^{\mathrm{vir}}=R\CH om_{\pi}(\CK,\CQ):=R \pi_* R\CH om(\CK,\CQ).
\end{equation}

\begin{lemma} \label{relquot2} There exists a natural identification
	\[\BT_{\mathrm{Quot}}^{\mathrm{vir}}\cong \CH om_{\pi}(\CK,\CQ) \oplus \CH om_{\pi}(\CK, \CQ)\mathbf t  [-1].\]
\end{lemma}
\textit{Proof.} The proof is similar to the one of Lemma \ref{lemma}. By pulling back via the projection $K_E \rightarrow E$, we extend $\CK$ and $\CQ$ from the zero section to the entire $K_E$. The extensions are denoted by  $\bar \CK$ and $\bar \CQ$. 
Consider the sequence on $K_E$, 
\[ 0 \rightarrow \CO(- E) \rightarrow  \CO_{K_{E}} \rightarrow \CO_{ E} \rightarrow 0,\]
which we tensor with $\bar \CK$ to obtain 
\[0 \rightarrow \bar \CK(- E) \rightarrow  \bar \CK \rightarrow \CK \rightarrow 0, \]
we then apply $R\CH om_{\pi}(-, \CQ)$ and take the associate long exact sequence to obtain
\begin{multline*}
	0 \rightarrow   \CH om_{\pi}(\CK, \CQ) \rightarrow   \CH om_{\pi}(\bar \CK, \CQ) \rightarrow  \CH om_{\pi}(\bar \CK(-E), \CQ)  \rightarrow \\
	\rightarrow  \CE xt^1_{\pi}(\CK, \CQ) \rightarrow\CE xt^1_ \pi( \bar \CK, \CQ)  \rightarrow \CE xt^1_\pi( \bar \CK(-E), \CQ )\rightarrow \dots 
\end{multline*}
Since $Q$ is scheme-theoretically supported on the zero section, the map  $\CH om_{\pi}(\bar \CK, \CQ) \rightarrow  \CH om_{\pi}(\bar \CK(-E), \CQ)$ is zero, hence the map $\CH om_{\pi}(\CK, \CQ) \rightarrow     \CH om_{\pi}(\bar \CK, \CQ)$ is an isomorphism. This also implies that $\CH om_{\pi}(\bar \CK(-E), \CQ)  \rightarrow\CE xt^1_{\pi}(\CK, \CQ)$ is injective.   Since we are considering quotients of a stable sheaf, $\mathrm{ext}^2(\CK,\CQ)=0$, hence 
\[\mathrm{ext}^0(K,Q)-\mathrm{ext}^1(K,Q)=\ch(K^\vee) \cdot \ch(Q)=0,\]
we conclude that  $\CH om_{\pi}(\bar \CK(-E), \CQ)  \rightarrow\CE xt^1_{\pi}(\CK, \CQ)$ is in fact an isomorphism. Using the $\BC^*_t$-equivariant identification 
\[\CO_{K_E}(-E)\cong \CO_{K_E}\mathbf t^{-1},\]
we obtain that 
\begin{equation*} 
	\CE xt^1_{\pi}(\CK, \CQ) \cong \CH om_{\pi}(\bar \CK(-E), \CQ)\cong \CH om_{\pi}(\bar \CK, \CQ)\mathbf t\cong \CH om_{\pi}(\CK, \CQ)\mathbf t. 
\end{equation*}
Finally, since $\mathrm{Quot}(\asf)$ is smooth (as scheme, it is a Quot scheme of quotients of a stable sheaf on a curve), we obtain that 
	\[\BT_{\mathrm{Quot}}^{\mathrm{vir}}\cong \CH om_{\pi}(\CK,\CQ) \oplus \CH om_{\pi}(\CK, \CQ)\mathbf t  [-1].\]
	\qed
\\

Now let $\CK^\rel$ and $\CQ^\rel$  be the universal relative kernel and relative quotient on $K_{C\times E} \times_C\mathrm{Quot}^\rel(\asf)$, and let $\pi \colon K_{C\times E} \times_C\mathrm{Quot}^\rel(\asf) \rightarrow \mathrm{Quot}^\rel(\asf)$ be the natural projection.   The $\BC_t^*$-equivariant relative obstruction theory of $\mathrm{Quot}^\rel(\asf)$ is given by the complex
\[ R \CH om_{\pi}(\CK^\rel,\CQ^\rel) . \] 
\begin{lemma} \label{relquot} There exists a natural identification
	\[\BT_{\mathrm{Quot^\rel}(\asf)}^{\mathrm{vir}}\cong \CH om_{\pi}(\CK^\rel,\CQ^\rel) \oplus \CH om_{\pi}(\CK^\rel,\CQ^\rel) \boxtimes\omega_C\mathbf t[-1].\]
\end{lemma}
 \textit{Proof.} As in Lemma \ref{relquot2} , there is a natural identification 
\[\CO_{K_{C\times E}}(-C\times E)_{|C\times E} \cong \omega^{\vee}_{C\times E} \mathbf t^{-1} \cong \omega^{\vee}_C \mathbf t^{-1},\]
and the sequence on $K_{C \times E}$,
\[ 0 \rightarrow \CO(-C\times E) \rightarrow  \CO_{K_{C\times E}} \rightarrow \CO_{C\times E} \rightarrow 0.\]
From which we obtain that 
\begin{multline*}
\CE xt^1_{\pi}(\CK^\rel,\CQ^\rel)\cong \CH om_{\pi}(\CK^\rel(-C\times E),\CQ^\rel)
 \cong \CH om_{\pi}(\CK^\rel,\CQ^\rel)  \boxtimes  \omega_C\mathbf t.
\end{multline*}
 So we conclude that the $C$-relative obstruction theory of $\mathrm{Quot}^\rel(\asf)$ can be expressed as follows,
\begin{equation*} 
	\BT_{\mathrm{Quot^\rel}}^{\mathrm{vir}}\cong \CH om_{\pi}(\CK^\rel,\CQ^\rel) \oplus \CH om_{\pi}(\CK^\rel,\CQ^\rel )\boxtimes\omega_C\mathbf t[-1]. 
\end{equation*}
\qed
\\

Let us compare two obstruction theories. With respect to the identification (\ref{ident1}), the sheaf  $\CH om_{\pi}(\CK^\rel,\CQ^\rel)$ is identified with the (pullback of) the sheaf $\CH om_{\pi}(\CK,\CQ)$,
\begin{equation} \label{tangent}
\CH om_{\pi}(\CK^\rel,\CQ^\rel) \cong \CH om_{\pi}(\CK,\CQ)\cong T_{\mathrm{Quot}(\asf)} ,
\end{equation}
 Hence the only difference between obstruction theories is the twist of the obstruction bundle by $\omega_C$ in Lemma \ref{relquot}. We conclude with the following corollary. 
 \begin{cor} \label{obst} \
\[ \Ob_{\mathrm{Quot}^\rel(\asf)}  \cong T_{\mathrm{Quot}(\asf)} \boxtimes \omega_C \mathbf t. \]
 	\end{cor}

\subsection{Group actions on Quot schemes}
The group $\Phi_\asf$ acts naturally on Quot schemes $\mathrm{Quot}(\asf )$. The stabilisers of the action are finite, as long as $\wsf\neq 0$. The quotient stack $[\mathrm{Quot}(\asf)/\Phi_\asf]$ is therefore a Deligne--Mumford stack.  Taking quotient respects  the identification \ref{ident1}, 
\[[\mathrm{Quot}^\rel(\asf)/\Phi_\asf]= [\mathrm{Quot}(\asf)/\Phi_\asf]\times C.\]
The obstruction theory of $[\mathrm{Quot}^\rel(\asf)/\Phi_\asf]$ is a descend of the obstruction theory of $[\mathrm{Quot}^\rel(\asf)]$.
By the preceding analysis of obstruction theories, the relative obstruction bundles of the two quotients are therefore  given by the descend of $\CH om_{\pi}(\CK,\CQ) \boxtimes  \omega_C\mathbf t$ and the descend of $\CH om_{\pi}(\CK,\CQ)\mathbf t$, respectively.  
 More precisely, let 
\[q \colon \mathrm{Quot}^\rel(\asf ) \rightarrow [\mathrm{Quot}^\rel(\asf )/\Phi_\asf] \]
be the quotient map, then  we have a $\Phi_\asf$-equivariant identification

\begin{equation*} \label{descend}
q^*\Ob_{[\mathrm{Quot}^\rel(\asf )/\Phi_\asf]}\cong\Ob_{\mathrm{Quot}^\rel(\asf )},
\end{equation*}
such that $\Ob_{\mathrm{Quot}^\rel(\asf)}$ is given the natural $\Phi_\asf$ equivariant structure. 
Hence by Corollary \ref{obst}, we obtain the following result. 
\begin{cor}\label{obst2} There is a $\Phi_\asf$-equivariant identification
\[\Ob_{[\mathrm{Quot}^\rel(\asf )/\Phi_\asf]} \cong T_{\mathrm{Quot}(\asf)}  \boxtimes  \omega_C\mathbf t . \] 
\end{cor}
We are now ready to determine the virtual degree. 
\begin{prop} \label{topch} For any $\rsf$ and $\asf$, we have
	\begin{equation*} \deg[\mathrm{Quot}^\rel(\asf )/\Phi_\asf]^{\mathrm{vir}}=(2g-2)e(T_{[\mathrm{Quot}(\asf )/\Phi_\asf]}).
		\end{equation*}
	\end{prop}
\textit{Proof.} 
Consider the map
\[ [\mathrm{Quot}(\asf )/\Phi_\asf] \rightarrow [\mathrm{pt}/\Phi_\asf],\] 
the associated sequence of tangent complexes is of the following form
\[ 0 \rightarrow  T_{\Phi_\asf} \rightarrow  T_{\mathrm{Quot}(\asf)} \rightarrow T_{[\mathrm{Quot}(\asf)/\Phi_\asf]} \rightarrow 0,\]
where by  $T_{\mathrm{Quot}(\asf)}$ we mean the descend of $T_{\mathrm{Quot}(\asf)}$ to $[\mathrm{Quot}(\asf)/\Phi_\asf]$. Since $T_{\Phi_\asf}$ is trivial, we obtain that 
\begin{align*}
	e(T_{[\mathrm{Quot}(\asf )}])&=0 \\
\mathrm{c}_{\rk-1}(T_{[\mathrm{Quot}(\asf )}])&= e(T_{[\mathrm{Quot}(\asf )/\Phi_\asf]}).
\end{align*}
Using Corollary \ref{obst2}, we therefore obtain
\begin{align*} \label{euler}
 [\mathrm{Quot}^\rel(\asf )/\Phi_\asf]^{\mathrm{vir}}&= e_{\BC^*_t}(\Ob_{[\mathrm{Quot}^\rel(\asf )/\Phi_\asf]}) \\
 &= e(T_{\mathrm{Quot}(\asf )})+ \mathrm{c}_{\rk-1}(T_{[\mathrm{Quot}(\asf )}])\cdot (\mathrm{c}_1(\omega_C)-t)+\dots \\
 &=  e(T_{[\mathrm{Quot}(\asf )/\Phi_\asf]})\cdot (\mathrm{c}_1(\omega_C)-t)+\dots.
 \end{align*}
Taking the degree, we arrive at the statement of the proposition, 
  \[ \deg[\mathrm{Quot}^\rel(\asf )/\Phi_\asf]^{\mathrm{vir}}= (2g-2)e(T_{[\mathrm{Quot}(\asf )/\Phi_\asf]}). \] 
  \qed 
  \\
  
   \noindent We now have to compute $e(T_{[\mathrm{Quot}(\asf )/\Phi_\asf]})$. The analysis will be split depending on the Chern cahracter of quotients.
\subsection{Relevant Quot schemes} \label{relevant}

For this section we assume that $\asf\neq0$. 
\subsubsection{} 
Assume 
\[u=(\rsf-1,\ksf).\]
Let
\[\mathsf{dim}:=\dim (\mathrm{Quot}(\asf, u ))= \rsf(\ksf-\asf)+\asf,\]
then $\mathrm{Quot}(\asf, u )$ is a $\p^{\mathsf{dim}-1}$-bundle over $\Jac(E)$ given by the natural projection 
\begin{align*}
	\mathrm{Quot}(\asf, u )& \rightarrow \Pic(E)\\
	[K \hookrightarrow G\twoheadrightarrow Q] &\mapsto  K.
\end{align*}
A fiber
\[ \p^{\mathsf{dim}-1} \hookrightarrow \mathrm{Quot}(\asf, u ) \]
is a slice of the $\Phi_\asf$-action on $\mathrm{Quot}(\asf, u )$. In other words, let $\Gamma_\ksf \subset \Phi_\asf$ be a finite subgroup that fixes $\p^{\mathsf{dim}-1}$, then we have the following diagram
\[
\begin{tikzcd}[row sep=scriptsize, column sep = scriptsize]
	& \p^{\mathsf{dim}-1}\arrow[d, "\pi_1"]  \arrow[r, hook] &  \mathrm{Quot}(\asf, u ) \arrow[d, "\pi_2"]\\
	& \left[\p^{\mathsf{dim}-1}/\Gamma_\ksf \right ] \arrow[r,"\cong"] &  \left[\mathrm{Quot}(\asf, u )/\Phi_\asf\right]
\end{tikzcd}
\]
The diagram above gives us that 
\begin{equation} \label{oddk}
 e(T_{[\mathrm{Quot}(\asf,u )/\Phi_\asf]})=\frac{e(\p^{\mathsf{dim}-1})}{|\Gamma_\ksf|}=\frac{\mathsf{dim}}{|\Gamma_\ksf|}.
 \end{equation}

It therefore remains to determine $\Gamma_\ksf$. The subgroup $\Gamma_\ksf \subset \Phi_\asf$ that fixes the slice is exactly the subgroup that fixes a line bundle $\CO({\asf-\ksf})$ of degree $\asf-\ksf$. A translation of $\CO(\asf-\ksf)$ by $\tau_p$ associated to a point $p \in E$ can be described as follows
\[\tau_p^*\CO(\asf-\ksf)=\CO(\asf-\ksf) \otimes L_p^{\asf-\ksf}, \]
where $L_p$ is line bundle corresponding to $p$ under the natural identification
\begin{align*}
	E&\xrightarrow{\cong} \Jac(E)\\
	p& \mapsto \CO(0_E) \otimes \CO(p)^{-1}=:L_p. 
\end{align*}
By the definition of $\Phi_\asf$ from (\ref{phi}), determining the stabiliser of $\CO(\asf-\ksf)$ in $\Phi_\asf$ therefore amounts to finding pairs 
\[(p, L^{\frac{\asf}{\rsf}}_p) \in E\times \Jac(E)\]
such that  
\begin{equation} \label{eq1}
	L_p^{\asf-\ksf} \otimes L_p^{-\frac{\asf}{\rsf}} =\CO_E.
\end{equation}
Raising the expression to the power of $\rsf$, we conclude that $L_p^\mathsf{dim}=\CO_E$.  As a group, $E$ is isomorphic to $\BR/\BZ \times \BR/\BZ$, hence with respect to the identification $E \cong \BR/\BZ \times \BR/\BZ$, we obtain that   
\[L_p=\left(\frac{n_1}{\mathsf{dim} }, \frac{n_2}{\mathsf{dim}} \right) \in \BR/\BZ \times \BR/\BZ,\]
$\rsf$-roots of $a_i=\frac{n_i}{\mathsf{dim}}$ are given by the following elements 
\[a^{\frac{1}{\rsf}}_i=\frac{n_i+h\mathsf{dim}}{\rsf\mathsf{dim}}, \quad h\in \{0,\dots,\rsf-1\}.\]
Only one satisfies the equation
\[\frac{(\asf-\ksf)n_i}{\mathsf{dim} }-\frac{\asf(n_i+h\rsf)}{\rsf\mathsf{dim} }=0 \in \BR/\BZ, \]
more specifically, $h$ is the uniquely defined by the following equation
\[\frac{ (n_i -\asf  h)\mathsf{dim}}{\rsf\mathsf{dim}}=0\in \BR/\BZ.\]
We therefore conclude that for all $p\in E$ there exists a unique root $L^{\frac{\asf}{\rsf}}_p$ which satisfies the equation $(\ref{eq1})$ and it must be $\mathsf{dim}$-torsion, hence $p$ is also $\mathsf{dim}$-torsion. We therefore obtain that  
\[ \Gamma_\ksf = E[\mathsf{dim}] \cong \BZ_{\mathsf{dim}}^{\oplus2},\]
in particular, 
\[ |\Gamma_\ksf|=\mathsf{dim}^2.\]
Proposition \ref{topch} and (\ref{oddk}) give us the following. 
\begin{cor} \label{lemma1} If $u=(\rsf-1,\ksf)$, then
	\begin{align*}
		&\deg([\mathrm{Quot}(\asf, u )^\rel/\Phi_\asf]^{\mathrm{vir}})= (2g-2)\dim(\mathrm{Quot}(\asf, u)^{-1}.
	\end{align*}
\end{cor}

\subsubsection{} \label{points}
Assume 

\[u=(0,\ksf).\]
In this case, 
\[\dim (\mathrm{Quot}(\asf, u ))= \rsf\ksf ,\]
$\mathrm{Quot}(\asf, u )$ admits a natural projection
\begin{align*}\mathrm{Quot}(\asf, u )& \rightarrow \Pic(E)\\
	[G \twoheadrightarrow Q ]& \mapsto \det(Q),
	\end{align*}
which provides a slice of the action of $\Phi_\asf$ on $\mathrm{Quot}(\asf, u )$. More specifically, let  
 \[ \mathrm{Quot}(\asf, u )_0 \hookrightarrow \mathrm{Quot}(\asf, u ) \]
be the fiber of the projection, let $\Gamma_\ksf \subset \Phi_\asf$ be the stabiliser of  $\mathrm{Quot}(\asf, u )_0$.  We obtain that
 
 \begin{equation} \label{evenk}
 	e(T_{[\mathrm{Quot}(\asf,u)/\Phi_1]})=\frac{e(\mathrm{Quot}(\asf, u )_0)}{|\Gamma_\ksf|}.
 \end{equation}
The subgroup $\Gamma_\ksf \in \Phi_\asf$ that fixes the slice is exactly the subgroup that fixes the determinant of a $0$-dimension sheaf of degree $\ksf$, this means that it consists of pairs 
\[(p, L^{\frac{\asf}{\rsf}}_p) \in E\times \Jac(E)\]
such that
\[ L_p^\ksf =\CO_E, \]
hence 
\[ |\Gamma_\ksf |= \rsf^2\ksf^2.\]

 Let us now determine the Euler characteristics of $\mathrm{Quot}(\asf, u )_0$. Firstly, any vector bundle on a curve can be deformed to a direct sum of line bundles. By the deformation invariance of the (virtual) Euler characteristics, we can assume that $G=\bigoplus^{i=\rsf}_{i=1} L_i $. In this case, $\mathrm{Quot}(\asf, u )_0$ admits a torus-action of $T= \prod^{i=\rsf }_{i=1}  \BC^*$ acting by scaling  line bundles. The associated fixed locus has the following description 
 \[ \mathrm{Quot}(\asf, u )_0^T= \coprod_{u_1\substack{+\dots +}u_\rsf=u}   \left( \prod^{i=\rsf}_{i=1} \mathrm{Quot}(L_i,u_i)\right)_0, \]
 we refer to \cite[Section 3]{MR} for more details in case of usual Quot schemes, which extend in a straightforward manner to our slices. Let us now analyse Euler characteristics of the components in the decomposition above. Firstly, if at least two class $u_k$ and $u_j$ are non-zero, then 
 $\left( \prod_i \mathrm{Quot}(L_i,u_i)\right)_0$ 
   admits an extra fixed-point-free action of $E$, therefore 
 \begin{align*}
 e\left( \left(\prod^{i=\rsf}_{i=1}\mathrm{Quot}(L_i,u_i)\right)_0\right)=0, \quad \text{ if }u_k\neq 0 \text{ and } u_j \neq 0 .
 \end{align*}
If only one class $u_j$ is non-zero, then
 \[\left( \prod^{i=\rsf}_{i=1} \mathrm{Quot}(L_i,u_i)\right)_0=  \mathrm{Quot}(L_{j},u_{j})_0= \p^{\ksf-1}.\]
  We therefore obtain 
  \begin{align*}
  e\left( \left( \prod^{i=\rsf}_{i=1}\mathrm{Quot}(L_i,u_i)\right)_0\right)=\ksf, \quad \text{ if only one } u_{j}\neq 0.
   	\end{align*}
Overall, 
 \begin{equation} \label{eulereven}
 	e(\mathrm{Quot}(\asf, u )_0)=e(\mathrm{Quot}(\asf, u )_0^T)=\rsf\ksf .
 	\end{equation}
Combining Proposition \ref{topch}, (\ref{evenk}) and (\ref{eulereven}), we obtain the following. 
\begin{cor} \label{lemma2} If $u=(0,\ksf)$, then
	\begin{align*} 
		\deg([\mathrm{Quot}(\asf, u )^\rel/\Phi_\asf]^{\mathrm{vir}})= (2g-2)\dim(\mathrm{Quot}(\asf, u)^{-1}.
	\end{align*}
\end{cor}
\bibliographystyle{amsalpha}
\bibliography{HiggsQM}

\providecommand{\bysame}{\leavevmode\hbox to3em{\hrulefill}\thinspace}
\providecommand{\MR}{\relax\ifhmode\unskip\space\fi MR }
\providecommand{\MRhref}[2]{%
  \href{http://www.ams.org/mathscinet-getitem?mr=#1}{#2}
}
\providecommand{\href}[2]{#2}
\begin{thebibliography}{{Zho}22}

\bibitem[CY22]{CY}
H.-D. Cao and S.-T. Yau (eds.), \emph{Differential geometry, {Calabi}-{Yau}
  theory, and general relativity. {Part} 2. {Lectures} and articles celebrating
  the 70th birthday of {Shing}-{Tung} {Yau}, {Harvard} {University},
  {Cambridge}, {MA}, {USA}, {May} 2019}, Surv. Differ. Geom., vol.~24,
  Somerville, MA: International Press, 2022.

\bibitem[Dij95]{Di}
R.~Dijkgraaf, \emph{Mirror symmetry and elliptic curves}, The moduli space of
  curves. Proceedings of the conference held on Texel Island, Netherlands
  during the last week of April 1994, Basel: Birkh{\"a}user, 1995,
  pp.~149--163.

\bibitem[DPS98]{DPS}
R.~Dijkgraaf, J.-S. Park, and B.~Schroers, \emph{{N=4 Supersymmetric Yang-Mills
  Theory on a K\"ahler Surface}}, arXiv:hep-th/9801066 (1998).

\bibitem[GK20]{GK2}
L.~G{\"o}ttsche and M.~Kool, \emph{Virtual refinements of the {Vafa}-{Witten}
  formula}, Commun. Math. Phys. \textbf{376} (2020), no.~1, 1--49.

\bibitem[GSY20]{GSY}
A.~Gholampour, A.~Sheshmani, and S.-T. Yau, \emph{Localized
  {Donaldson}-{Thomas} theory of surfaces}, Am. J. Math. \textbf{142} (2020),
  no.~2.

\bibitem[HT03]{HT}
T.~Hausel and M.~Thaddeus, \emph{Mirror symmetry, {Langlands} duality, and the
  {Hitchin} system}, Invent. Math. \textbf{153} (2003), no.~1, 197--229.

\bibitem[Joy21]{J}
D.~Joyce, \emph{{Enumerative invariants and wall-crossing formulae in abelian
  categories}}, arXiv:2111.04694 (2021).

\bibitem[Laa20]{Laa}
Ties Laarakker, \emph{Monopole contributions to refined {Vafa}-{Witten}
  invariants}, Geom. Topol. \textbf{24} (2020), no.~6, 2781--2828.

\bibitem[LW23]{LeeW}
S.~Lee and Y.~Wen, \emph{{Quasisections to moduli spaces of sheaves}}, in
  preparation.

\bibitem[MM21]{MM}
J.~Manschot and G.~W. Moore, \emph{{Topological correlators of $SU(2)$,
  $\mathcal{N}=2^*$ SYM on four-manifolds}}, arXiv:2104.06492 (2021).

\bibitem[Moc09]{Moch}
T.~Mochizuki, \emph{Donaldson type invariants for algebraic surfaces.
  {Transition} of moduli stacks}, Lect. Notes Math., vol. 1972, Berlin:
  Springer, 2009.

\bibitem[MR22]{MR}
S.~Monavari and A.~T. Ricolfi, \emph{On the motive of the nested {Quot} scheme
  of points on a curve}, J. Algebra \textbf{610} (2022), 99--118.

\bibitem[Nes21]{N}
D.~Nesterov, \emph{Quasimaps to moduli spaces of sheaves}, arXiv:2111.11417
  (2021).

\bibitem[Nes23]{NHiggs}
\bysame, \emph{{Enumerative mirror symmetry for moduli spaces of Higgs bundles
  and S-duality}}, arXiv:2302.08379 (2023).

\bibitem[Obe21]{Ob}
G.~Oberdieck, \emph{{Multiple cover formulas for K3 geometries, wallcrossing,
  and Quot schemes}}, arXiv:2111.11239 (2021).

\bibitem[Oko19]{Ok}
Andrei Okounkov, \emph{Takagi lectures on {Donaldson}-{Thomas} theory}, Jpn. J.
  Math. (3) \textbf{14} (2019), no.~1, 67--133.

\bibitem[PT16]{KKV}
R.~Pandharipande and R.~P. Thomas, \emph{The {Katz}-{Klemm}-{Vafa} conjecture
  for {{\(K3\)}} surfaces}, Forum Math. Pi \textbf{4} (2016), 111.

\bibitem[Tho20]{Th20}
R.~P. Thomas, \emph{Equivariant {{\(K\)}}-theory and refined {Vafa}-{Witten}
  invariants}, Commun. Math. Phys. \textbf{378} (2020), no.~2, 1451--1500.

\bibitem[TT20]{TT}
Y.~Tanaka and R.~P. Thomas, \emph{Vafa-{Witten} invariants for projective
  surfaces. {I}: {Stable} case}, J. Algebr. Geom. \textbf{29} (2020), no.~4,
  603--668.

\bibitem[{Zho}22]{YZ}
Y.~{Zhou}, \emph{{Quasimap wall-crossing for GIT quotients}}, {Invent. Math.}
  \textbf{227} (2022), no.~2, 581--660.

\end{thebibliography}

	\end{document}